\documentclass[11pt]{amsart}
\usepackage{amsmath, amssymb, comment, epsfig, amsthm, psfrag}

\newcommand{\K}{\mathcal{K}}
\newcommand{\Gr}{\mathcal{G}r}
\newcommand{\wl}{X^*}
\newcommand{\cwl}{X_*}
\newcommand{\D}{\mathrm{D}}
\newcommand{\C}{\mathbb{C}}
\newcommand{\Z}{\mathbb{Z}}
\newcommand{\R}{\mathbb{R}}
\newcommand{\N}{\mathbb{N}}
\newcommand{\val}{\operatorname{val}}
\newcommand{\dual}{\star}
\newcommand{\kpf}{\operatorname{kpf}}
\newcommand{\row}{\Gamma}
\newcommand{\conv}{\operatorname{conv}}
\newcommand{\bigdot}{\bullet}
\newcommand{\wi}{\mathbf{i}}
\newcommand{\wj}{\mathbf{j}}
\newcommand{\realt}{\mathfrak{t}_\R}
\newcommand{\fund}{\Lambda}
\newcommand{\Hom}{\operatorname{Hom}}
\newcommand{\Verti}{\operatorname{Vert}}
\newcommand{\rdim}{\operatorname{rdim}}
\newcommand{\Grl}{\mathcal{G}rl}
\newcommand{\Span}{\operatorname{span}}
\newcommand{\oo}{\mathcal{O}}

\newcommand{\Comp}{\operatorname{Comp}}

\newtheorem{Theorem}{Theorem}[section]
\newtheorem{Specialthm}{Theorem}
\newtheorem{Proposition}[Theorem]{Proposition}
\newtheorem{Lemma}[Theorem]{Lemma}
\newtheorem{Corollary}[Theorem]{Corollary}

\newtheorem{Question}{Question}

\theoremstyle{definition}
\newtheorem{Example}{Example}

\begin{document}

\title{Mirkovi\'c-Vilonen cycles and polytopes}
\author{Joel Kamnitzer}

\begin{abstract}
We give an explicit description of the Mirkovi\'c-Vilonen cycles on the affine Grassmannian for arbitrary complex reductive groups.  We also give a combinatorial characterization of the MV polytopes.  We prove that a polytope is an MV polytope if and only if it a lattice polytope whose defining hyperplanes are parallel to those of the Weyl polytopes and whose 2-faces are rank 2 MV polytopes.  As an application, we give a bijection between Lusztig's canonical basis and the set of MV polytopes.  
\end{abstract}

\maketitle

\tableofcontents

\section{Introduction}
\subsection{Background}
Let $ G $ be a complex connected reductive group and let $ G^\vee $ be its Langlands dual group.  Let $ \K = \mathbb{C}((t)) $ denote the field of Laurent series and let $ \mathcal{O} = \mathbb{C}[[t]] $ denote the ring of power series.  The quotient $ \Gr = G(\K) / \, G(\mathcal{O}) $ is called the \textbf{affine Grassmannian}.  The geometric Satake correspondence of Lusztig \cite{L}, Ginzburg \cite{G}, Beilinson-Drinfeld \cite{BD}, and Mirkovi\'c-Vilonen \cite{MV2} provides a connection between the geometry of $ \Gr $ and the representation theory of $ G^\vee $.

\begin{Specialthm}[Lusztig]
For each $ \lambda \in \cwl^+ $, the set of dominant weights of $ G^\vee$, there exists a subvariety $ \Gr^\lambda $ of $ \Gr $ such that $ \mathrm{IH}(\Gr^\lambda) \cong V_\lambda $. 
\end{Specialthm}

Here $ \mathrm{IH}(\Gr^\lambda) $ denotes the intersection homology of $ \Gr^\lambda $ and $V_\lambda $ denotes the irreducible representation of $ G^\vee $ of highest weight $\lambda $. 

As a simple example, take $ G = GL_n $ (in this case $ G^\vee = GL_n $) and $ \lambda = (1, \dots, 1, 0, \dots ,0) $ (where there are $ k $ 1s).  Then $ \Gr^\lambda \cong \mathrm{Gr}(k,n) $, the usual Grassmannian of $ k $-planes in $ \C^n $.  Since $ \mathrm{Gr}(k,n) $ is smooth, $ \mathrm{IH}(\Gr^\lambda) = \mathrm{H}(\mathrm{Gr}(k,n)) $.  Recall that the homology of $ \mathrm{Gr}(k,n) $ has a basis given by the Schubert varieties, which are naturally indexed by $ k $ element subsets of $ \{1, \dots, n \} $.

In this case, the right hand side $ V_\lambda$ is the representation of $ GL_n $ on $ \Lambda^k \C^n $, which also has a basis indexed by $ k $-element subsets of $ \{1, \dots, n\}$.  The geometric Satake correspondence says that this is not a coincidence, but rather part of a larger pattern which holds for all finite dimensional representations of complex reductive groups.

Mirkovi\'c-Vilonen extended Lusztig's work as follows.
\begin{Specialthm}[Mirkovi\'c-Vilonen]
There exists a family of subvarieties of the affine Grassmannian, called Mirkovi\'c-Vilonen cycles, such that the subset lying in $ \Gr^\lambda $ forms a basis for $ \mathrm{IH}(\Gr^\lambda) $.
\end{Specialthm}
Hence we get a basis for $ V_\lambda $ indexed by MV cycles.  In the above example, the MV cycles are exactly the Schubert varieties.

These theorems motivate the following question:
\begin{Question} \label{qu:comb}
Can we use the MV cycles in $ \Gr^\lambda $ to understand the combinatorics of bases for the representation $ V_\lambda $?
\end{Question}
In our simple example, we can use the Schubert varieties in $ \mathrm{Gr}(k,n) $ to see that $ \Lambda^k \C^n $ has a basis indexed by the $ k $ element subsets of $ \{1, \dots, n\}$. 

Some attempts have been made to give a combinatorial description of the MV cycles.  The problem is that the MV cycles are mysterious, since they are defined as the components of intersections of opposite ``semi-infinite orbits''.  Gaussent-Littelmann \cite{GL} associated an MV cycle to each Littelmann path, by considering certain resolutions of $ \Gr^\lambda $.  

A different approach is due to Anderson \cite{jared2}.  He proposed understanding MV cycles by looking at their moment polytopes, which he called MV polytopes.  Anderson used the above results of Lusztig and Mirkovi\'c-Vilonen to show that MV polytopes could be used to count weight and tensor product multiplicities for $ G^\vee $.  However, he could not give a characterization of the MV polytopes since he did not have an explicit description of the MV cycles. 

Anderson-Kogan \cite{jaredmisha} studied MV cycles for $ GL_n $ by means of the lattice model for $ \Gr $.  They gave a recipe for producing MV cycles and polytopes for $ GL_n $, but not an explicit description of the cycles and polytopes.  

\subsection{Main result}
In this paper, we give an explicit combinatorial description of the MV cycles and polytopes uniform across all types.  We begin with the notions of ``pseudo-Weyl polytope'' and ``GGMS stratum'' (see sections \ref{se:pWeyl}, \ref{se:GGMS}).  A pseudo-Weyl polytope is a lattice polytope whose defining hyperplanes are parallel to those of the Weyl polytopes.  A GGMS stratum, whose moment map image is a pseudo-Weyl polytope, is the intersection of semi-infinite cells, one for each element of the Weyl group.  A pseudo-Weyl polytope and a GGMS stratum are each described by a collection of integers, one for each ``chamber weight''.  On the polytope side, these integers give the positions of the defining hyperplanes, while on the GGMS stratum side, they are the values of certain constructible functions (section \ref{se:Dgam}).  More concretely, if $ \big( M_\gamma \big)_{\gamma \in \Gamma} $ is such a collection of integers, then
\begin{equation*}
P(M_\bigdot) := \{ \alpha \in \realt : \langle \alpha, \gamma \rangle \ge M_\gamma \text{ for all } \gamma \} \quad A(M_\bigdot) := \{ L \in \Gr : \D_\gamma(L) = M_\gamma \text{ for all } \gamma \}
\end{equation*}
are the corresponding pseudo-Weyl polytope and GGMS stratum.

The important point is to determine for which collections of integers is the closure of the resulting GGMS stratum an MV cycle.  The key idea is that our constructible functions are closely related to the valuations of the generalized minors of Berenstein-Zelevinsky \cite{BZschub} and that the Pl\"ucker relations hold among these generalized minors.  Thus we are lead to the tropical form of these relations, which is obtained by replacing $ + $ with $ \min $ and $ \times $ with $ + $ in these relations (see section \ref{se:tpr}).  

\begin{Specialthm}[Theorem \ref{th:BZcycle}]
If $ M_\bigdot $ satisfies the ``tropical Pl\"ucker relations'' and certain ``edge inequalities'', then $ \overline{A(M_\bigdot)} $ is an MV cycle and $ P(M_\bigdot) $ is an MV polytope.  Moreover all MV cycles and polytopes arise this way.  
\end{Specialthm}

The following corollary follows from the form of the ``tropical Pl\"ucker relations''.
\begin{Specialthm}
A pseudo-Weyl polytope is an MV polytope if and only if every 2-face is an MV polytope of the appropriate rank 2 group.  The MV polytopes for $ SL_3 $ and $ Sp_4 $ are given in figures \ref{fig:sl3} and \ref{fig:sp4}.
\end{Specialthm}

Section \ref{se:decomp} and \ref{se:overlap} are devoted to the proof of Theorem \ref{th:BZcycle}.   In section \ref{se:decomp}, we explain how each reduced word $ \wi $ for $ w_0 $ gives a decomposition of the affine Grassmannian into irreducible pieces according to $\wi$-Lusztig datum.  We prove (Theorem \ref{th:Ldirr}) that the closures of these pieces are the MV cycles.  In this section, we use the results of Berenstein-Fomin-Zelevinsky \cite{BFZA,BZschub,FZ} concerning generalized minors.
In section \ref{se:overlap}, we consider the overlap of decompositions for different $\wi$.  The key is to first consider reduced words $ \wi, \wi'$ which differ by a braid move (section \ref{se:loc}).  Here we use a result of Lusztig and Berenstein-Zelevinsky on the comparison between different parametrizations of the upper triangular subgroup of $ G $.  Using this knowledge, we are able to prove that the MV cycles are as described in Theorem \ref{th:BZcycle}.

\subsection{Applications}
After proving this main theorem, we give a number of applications.  First we consider the problem of decomposing MV polytopes under Minkowski sum (section \ref{se:mink}).  In low rank cases, Anderson \cite{jared2} gave certain ``prime'' MV polytopes which he conjectured generated all the MV polytopes under Minkowski sums.  We show that for any group $ G $, there exists such a finite set of prime MV polytopes and moreover we show how to find these prime polytopes (Theorem \ref{th:MVprime}).

Combining our result with the work of Lusztig \cite{Lbook} and Berenstein-Zelevinsky \cite{BZtpm}, shows that there is a bijection between Lusztig's canonical basis and the set of MV polytopes (Theorem \ref{th:can}).  In the case of MV cycles, the tropical Pl\"ucker relations appear naturally (see section \ref{se:tpr}), whereas their appearance in \cite{BZtpm} to describe the canonical basis was unexpected.  Thus, we have bijections
\begin{equation} \label{eq:bij}
\mathcal{B} \longleftrightarrow \mathcal{P} \longleftrightarrow \mathcal{M}
\end{equation}
where $ \mathcal{B} $ denotes the canonical basis, $ \mathcal{P} $ denotes the set of MV polytopes, and $ \mathcal{M} $ denotes the set of MV cycles.  In \cite{MVcrystal}, we show that these bijections are isomorphisms of crystals with respect to the Kashiwara-Lustzig crystal structure on the canonical basis and the Braverman-Finkelberg-Gaitsgory crystal structure on the set of MV cycles.

Another important application of our main result is to answer Question \ref{qu:comb}.  Using the work of Mirkovic-Vilonen \cite{MV} and Anderson \cite{jared2}, we give a combinatorial description of the BZ data which index the MV basis for $ V_\lambda $ (Theorem \ref{th:MVwtm}).  In \cite{BZtpm}, Berenstein-Zelevinsky gave the BZ data which index the canonical basis for $ V_\lambda$.  These two sets are the same, even though there is a subtle difference in their descriptions.  Finally, we use the work of Anderson to give a tensor product multiplicity formula in terms of counting BZ data.

There is a close connection between our work and the Anderson-Kogan description of MV cycles and polytopes for $ GL_n $.  In fact, their work served as an important source of motivation.  The details of this connection are explained in section \ref{se:typeA}.  In particular, we show that their methods of producing MV cycles and polytopes from Kostant pictures fits into our framework (Theorems \ref{th:AKbij} and \ref{th:mycollapse}).

\subsection{Acknowledgements}
I would first like to thank my advisor Allen Knutson.  His encouragement and suggestions have proved valuable at many key stages of this project.

I thank Alexander Braverman and Peter Littelmann for lectures and conversations which led me to start thinking about MV cycles.  I am grateful to David Speyer for sharing an idea that was crucial to the beginning of this work.  I also thank Jared Anderson for beginning the study of MV polytopes and for many conversations and ideas.   In thinking about MV cycles and polytopes, I also benefited from conversations with David Ben-Zvi, Arkady Berenstein, Edward Frenkel, Tom Graber, Andr\'e Henriques, Misha Kogan, Ivan Mirkovi\'c, Scott Morrison, David Nadler, Arun Ram, Bernd Sturmfels, Peter Tingley, Kari Vilonen, Soroosh Yazdani, and Andrei Zelevinsky.  

I am especially grateful to Allen Knutson and Peter Tingley for their careful reading of this text.

During this work, I was supported financially by an NSERC postgraduate scholarship.

\section{Main definitions} \label{se:MVcyc}

\subsection{Notation}
If $ G $ is complex connected reductive group, then its affine Grassmannian is the disjoint union of $ \pi_1(G) $ many copies of the affine Grassmannian of the simply-connected semisimple group with the same root system as $ G $.  So here we only consider the case $ G $ connected simply-connected semisimple.  As another simplification, we consider only the case of $ G $ singly and doubly-laced.  Extending our results to include $ G_2 $ factors is quite simple; it just requires including the extra cases of $ a_{ij} = -3 $ and $ a_{ji} = -3 $ in the statement of the tropical Pl\"ucker relations (section \ref{se:tpr}) and in Propositions \ref{th:paramtrans} and \ref{th:valtrans}.  The case $a_{ij} = -3 $ appears in \cite{BZschub} and the case $ a_{ji} = -3 $ can be easily derived from there.

Let $ G $ be a connected simply-connected semisimple complex group. 

Let $ T $ be a maximal torus of $ G $ and let $ \wl = \Hom(T, \C^\times), \cwl = \Hom(\C^\times, T) $ denote the weight and coweight lattices of $ T $.  Let $ \Delta \subset \wl $ denote the set of roots of $ G $.  Let $ W = N(T)/ T $ denote the Weyl group.

Let $ B $ be a Borel subgroup of $ G $ containing $ T $.  Let $ \alpha_1, \dots, \alpha_r $ and $ \alpha_1^\vee, \dots, \alpha_r^\vee $ denote the simple roots and coroots of $ G $ with respect to $ B $.  Let $ N $ denote the unipotent radical of $ B $.  Let $ \fund_1, \dots, \fund_r $ be the fundamental weights.  Let $ I = \{ 1, \dots, r \} $ denote the vertices of the Dynkin diagram of $ G $.  Let $ a_{ij} = \langle \alpha_i^\vee, \alpha_j \rangle $ denote the Cartan matrix.  Let $ \rho := \sum \fund_i, \ \rho^\vee := \sum \fund_i^\vee $ be the Weyl and dual Weyl vectors.

Let $ s_1, \dots, s_r \in W $ denote the simple reflections.  Let $ e $ denote the identity in $ W $ and let $ w_0 $ denote the longest element of $ W $.  Let $ m $ denote the length of $ w_0 $ or equivalently the number of positive roots.  We will also need the Bruhat order on $ W $, which we denote by $ \ge $. 

We also use $ \ge $ for the usual partial order on $ \cwl $, so that $ \mu \ge \nu $ if and only if $ \mu - \nu $ is a sum of positive coroots.  More generally, we have the twisted partial order $ \ge_w $, where $ \mu \ge_w \nu $ if and only if $ w^{-1}  \cdot \mu \ge w^{-1} \cdot \nu $.

Let $ \realt := \cwl \otimes \R $ (the Lie algebra of the compact form of $ T $).  For each $ w$, we extend $ \ge_w $ to a partial order on $ \realt $, so that $ \beta \ge_w \alpha $ if and only if $ \langle \beta - \alpha, w \cdot \fund_i \rangle \ge 0 $ for all $ i $.

For each $ i \in I $, let $ \psi_i : SL_2 \rightarrow G $ be denote the $ i$th root subgroup of $ G $.

For $ w \in W $, let $ \overline{w} $ denote the lift of $ w $ to $ G $, defined using the lift of $ \overline{s_i} := \psi_i \Big( \big[ \begin{smallmatrix} 0 & 1 \\ -1 & 0 \\ \end{smallmatrix} \big] \Big) $.

A \textbf{reduced word} for an element $ w \in W $ is a sequence of indices $ \textbf{i} = (i_1, \dots, i_k) \in I^k$ such that $ w = s_{i_1} \cdots s_{i_k} $ is a reduced expression.

Let $ \kpf $ denote the \textbf{Kostant partition function} on $ \cwl $, so that $ \kpf(\mu) $ is the number of ways to write $ \mu $ as a sum of positive coroots.

If $ X $ is any variety, we write $ \Comp(X) $ for the set of components of $ X $.

\subsection{Affine Grassmannian}
For the purposes of this paper, it will be convenient to write the affine Grassmannian as the left quotient $ \Gr = G(\mathcal{O}) \setminus G(\K) $.  We view $ \Gr $ as an ind-scheme over $\C$ whose set of $ \C $ points is $ G(\mathcal{O}) \setminus G(\K) $.  Similarly, we view $ G(\K), N(\K), \K^m $ as ind-schemes over $ \C $.  More explicitly, they are the results of applying the formal loop space functor to $ G, N, \C^m $.  For more details, see \cite[Sections 11.3.3, 20.3.3]{BF}.

A coweight $ \mu \in \cwl $ gives a homomorphism $ \mathbb{C}^\times \rightarrow T $ and hence an element of $ \Gr$.  We denote the corresponding element $ t^\mu $.  It is easy to see that these $ t^\mu $ are the fixed points for the action of $ T(\mathbb{C}) $ on $\Gr$. 

For $ w \in W $, let $ N_w = w N w^{-1} $.  For $ w \in W $ and $ \mu \in \cwl $ define the \textbf{semi-infinite cells}
\begin{equation}
 S_w^\mu := t^\mu N_w(\K). 
\end{equation}
To a certain extent, these semi-infinite cells behave like the Schubert cells on a finite dimensional flag variety.  In particular, they are each attracting cells for a certain $\C^\times$ action on $\Gr $.  The choice of $ w \in W $ gives us a map $ w \cdot \rho^\vee : \C^\times \rightarrow T $ and we have
\begin{equation} \label{eq:Sattr}
S^\mu_w = \{ L \in \Gr: \lim_{s \rightarrow \infty} L \cdot (w \cdot \rho^\vee)(s) = t^\mu \}.
\end{equation} 

The semi-infinite cells have the simple containment relation (see \cite{MV})
\begin{equation} \label{eq:Scon}
\overline{S_w^\mu} = \bigcup_{\nu \ge_w \mu} S_w^\nu.
\end{equation}

\begin{Lemma} \label{th:Smeet}
If $S_w^\mu \cap S_v^\nu \ne \emptyset $ then $ \nu \ge_w \mu$.
\end{Lemma}

\begin{proof}
Let $ L \in S_w^\mu \cap S_v^\nu $.  Then by (\ref{eq:Sattr}), $ t^\nu  = \lim_{s \rightarrow \infty} L \cdot (v \cdot \rho^\vee)(s) $.  Since $ S^\mu_w $ is $ T $-invariant, this shows that $ t^\nu \in \overline{S^\mu_w} $.  So by (\ref{eq:Scon}), $ \nu \ge_w \mu $.
\end{proof}

Let $ \mu_1, \mu_2 $ be coweights with $ \mu_1 \le \mu_2 $.  Following Anderson \cite{jared2}, a component of $ \overline{S_e^{\mu_1} \cap S_{w_0}^{\mu_2}} $ is called an \textbf{MV cycle} of coweight $(\mu_1, \mu_2) $.  It is well-known that this intersection is finite dimensional.  (In fact, it is known that this intersection has pure dimension $ \langle \mu_2 - \mu_1, \rho \rangle $, but we will not need this fact.)

Note that $ \cwl $ acts on $ \Gr $ by $ \nu \cdot L := L \cdot t^\nu $.  Since $ T $ normalizes $ N_w $, we see that $ \nu \cdot S_w^\mu = S^{\mu+\nu}_w $.  So if $ A $ is a component of $ \overline{S_e^{\mu_1} \cap S_{w_0}^{\mu_2}} $, then $ \nu \cdot A $ is a component of $ \overline{S_e^{\mu_1 + \nu} \cap S_e^{\mu_2 + \nu}} $.  So $ \cwl $ acts on the set of all MV cycles. The orbit of an MV cycle of coweight $ (\mu_1, \mu_2) $ is called a \textbf{stable MV cycle} of coweight $\mu_2 - \mu_1 $.  Note that a stable MV cycle of coweight $ \mu $ has a unique representative of coweight $(\nu, \nu + \mu) $ for any coweight $\nu $.

Let $ \mathcal{M} $ denote the set of stable MV cycles and let $ \mathcal{M}(\mu) $ denote the set of those of coweight $\mu $.  It is well-known that there are $ \kpf(\mu) $ stable MV cycles of coweight $ \mu $ (for example this follows from \cite[section 13]{BFG}, or from \cite{jared2}).

Following Anderson \cite{jared2}, given a $ T $-invariant closed subvariety $ A $ of the affine Grassmannian, let $ \Phi(A) \subset \realt $ be the convex hull of $ \{ \mu \in \cwl : t^\mu \in A \}$.  By \cite{jared2}, this is the moment polytope for the $ T $ action on the affine Grassmannian.  

For example, by (\ref{eq:Scon}), we see that $ \Phi(\overline{S_w^\mu}) =  C_w^\mu := \{ \alpha \in \realt : \alpha \ge_w \mu \} = \{ \alpha : \langle \alpha, w \cdot \fund_i \rangle \ge \langle \mu, w \cdot \fund_i \rangle $ for all $ i \} $.

If $ A $ is an MV cycle of coweight $ (\mu_1, \mu_2) $, then we say that $ \Phi(A) $ is an \textbf{MV polytope of coweight} $ (\mu_1, \mu_2) $.  The action of $ \cwl $ on the set of MV cycles gives an action of $ \cwl $ on the set of MV polytopes.  In fact, it is easy to see that $ \nu \cdot P = P + \nu $.  The orbit of an MV polytope of coweight $ (\mu_1, \mu_2) $ is called a \textbf{stable MV polytope} of coweight $ \mu_2 - \mu_1 $.  Let $ \mathcal{P} $ denote the set of stable MV polytopes.

\subsection{Pseudo-Weyl polytopes} \label{se:pWeyl}
We will start our investigation by examining a larger family of polytopes, called pseudo-Weyl polytopes.  We will show how to pick out the MV polytopes from the pseudo-Weyl polytopes.  The idea that all MV polytopes should be pseudo-Weyl polytopes is due to Anderson.

For $ \lambda \in \cwl^+ $, $ W_\lambda = \conv (W \cdot \lambda) \subset \realt $ is called the $ \lambda$-\textbf{Weyl polytope}.  Recall that the Weyl polytope $ W_\lambda $ can be described in three different ways.  It is the convex hull of the orbit of $ \lambda $, it is the intersection of translated and reflected cones, and it is the intersection of half spaces.  In particular, 
\begin{equation*}
W_\lambda = \bigcap_w C_w^{w \cdot \lambda} = \{ \alpha \in \realt : \langle \alpha, w \cdot \fund_i \rangle \ge \langle w_0 \cdot \lambda, \fund_i \rangle \text{ for all $w \in W $ and $ i \in I$} \}.
\end{equation*}
Following Berenstein-Zelevinsky \cite{BZschub}, we call a weight $ w \cdot \fund_i $ a \textbf{chamber weight} of level $ i $.  So the chamber weights $ \Gamma := \bigcup_{w \in W, i \in I} w \cdot \fund_i $ are dual to the hyperplanes defining any Weyl polytope.

Suppose we are given a collection of coweights $ \mu_\bigdot = (\mu_w)_{w \in W} $ such that $ \mu_v \ge_w \mu_w $ for all $ v,w \in W $.  Then we can form the polytope 
\begin{equation*}
P(\mu_\bigdot) := \bigcap_w C_w^{\mu_w}
\end{equation*}
A \textbf{pseudo-Weyl polytope} is any polytope of this form.  

Pseudo-Weyl polytopes also admit a description in terms of intersecting half spaces.

Let $ M_\bigdot = \big( M_\gamma \big)_{\gamma \in \Gamma} $ be a collection of integers, one for each chamber weight.  Given such a collection, we can form $ P(M_\bigdot) := \{ \alpha \in \realt : \langle \alpha, \gamma \rangle \ge M_\gamma \text{ for all } \gamma \in \Gamma \} $.  This is the polytope made by translating the hyperplanes defining the Weyl polytopes to distances $ M_\gamma $ from the origin.

\begin{Proposition} \label{th:BZtoGGMS}
Let $ \mu_\bigdot = \big( \mu_w \big)_{w \in W} $ be a collection of coweights such that $ \mu_v \ge_w \mu_w $ for all $ v, w $.  Then the set of vertices of $ P(\mu_\bigdot) $ is the collection $ \mu_\bigdot $ (which may have repetition).

A pseudo-Weyl polytope has defining hyperplanes dual to the chamber weights.  In particular, if $ P $ is a pseudo-Weyl polytope with vertices $ \mu_\bigdot $, then $P = P(M_\bigdot) $ where
\begin{equation} \label{eq:BZtoGGMS}
 M_{w \cdot \fund_i } = \langle \mu_w, w \cdot \fund_i \rangle.
\end{equation}

Moreover, the $ M_\bigdot $ satisfy the following condition which we call the \textbf{edge inequalities}.
For each $ w \in W $ and $ i \in I $, we have:
\begin{equation} \label{eq:nondeg}
M_{ws_i \cdot \fund_i} + M_{w \cdot \fund_i} + \sum_{j \ne i} a_{j i} M_{w \cdot \fund_j} \le 0
\end{equation}

Conversely, suppose that a collection of integers $ \big(M_\gamma \big)_{\gamma\in \Gamma } $ satisfies the edge inequalities.  Then the polytope $ P(M_\bigdot) $ is pseudo-Weyl polytope with vertices given by 
\begin{equation} \label{eq:mufromM}
\mu_w = \sum_i M_{w \cdot \fund_i} w \cdot \alpha^\vee_i.
\end{equation}

\end{Proposition}

\begin{figure}
\begin{center}
\psfrag{me}{$\mu_e$} \psfrag{ms1}{$\mu_{s_1}$} \psfrag{ms1s2}{$\mu_{s_1 s_2}$}
\psfrag{ms1s2s1}{$\mu_{s_1s_2s_1}$} \psfrag{ms2s1}{$\mu_{s_2s_1}$} \psfrag{ms2}{$\mu_{s_2}$}
\psfrag{f1}{$\scriptstyle{\fund_1}$} \psfrag{s1f1}{$\scriptstyle{s_1 \cdot \fund_1}$} \psfrag{s1s2f2}{$\scriptstyle{s_1 s_2 \cdot \fund_2}$}
\psfrag{f2}{$\scriptstyle{\fund_2}$} \psfrag{s2f2}{$\scriptstyle{s_2 \cdot \fund_2}$} \psfrag{s2s1f1}{$\scriptstyle{s_2 s_1 \cdot \fund_1}$}
\epsfig{file=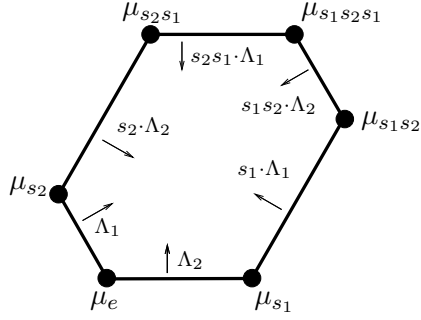, height=4cm}
\caption{A pseudo-Weyl polytope for $SL_3 $.}
\label{fig:pWeylsl3}
\end{center}
\end{figure}

 Figure \ref{fig:pWeylsl3} shows an example of a pseudo-Weyl polytope for $ G = SL_3 $ with vertices and chamber weights labelled.

In appendix \ref{se:appendix}, we will introduce the notion of dual fan of a polytope and will show that pseudo-Weyl polytopes are exactly those lattice polytopes whose dual fan is a coarsening of the Weyl fan.  We will also give a proof of Proposition \ref{th:BZtoGGMS}. 

Let $ P $ be a pseudo-Weyl polytope, $ P = \conv(\mu_\bigdot) = P(M_\bigdot) $.  For any $ w \in W , i \in I $, there is an edge connecting $ \mu_w $ and $ \mu_{ws_i}$.  We have
\begin{equation} \label{eq:length}
\mu_{w s_i} - \mu_w = c \, w \cdot \alpha^\vee_i, \ \text{ where } c = - M_{w \cdot \fund_i} - M_{w s_i \cdot \fund_i} - \sum_{j\ne i} a_{ji} M_{w \cdot \fund_j}.
\end{equation}
We call $ c $ the \textbf{length} of the edge from $ \mu_w $ to $ \mu_{w s_i} $.  Note that it is the negative of the left hand side of (\ref{eq:nondeg}).

\subsection{GGMS strata} \label{se:GGMS}
The geometric version of the pseudo-Weyl polytopes are the Gelfand-Goresky-MacPherson-Serganova (GGMS) strata on the affine Grassmannian.  These GGMS strata will be our candidates to be MV cycles.  These GGMS strata on the affine Grassmannian were first considered as potential MV cycles by Anderson-Kogan \cite{jaredmisha}.

Given any collection $ \mu_\bigdot = \big( \mu_w \big)_{w \in W} $ of coweights, we can form the \textbf{GGMS stratum}
\begin{equation} \label{eq:GGMSstratum}
A(\mu_\bigdot) := \bigcap_{w \in W} S_w^{\mu_w}.
\end{equation}

By Lemma \ref{th:Smeet}, this intersection is empty unless $ \mu_v \ge_w \mu_w $ for all $ v, w $.  So we will only consider such collections.

We will prove that each MV cycle is the closure of $ A(\mu_\bigdot) $ for an appropriate $ \mu_\bigdot $.   Once we know which of these are MV cycles, we will also know the MV polytopes, since we have the following easy lemma, which is a version of Theorem D from \cite{jaredmisha}.

\begin{Lemma} \label{th:momentpoly}
Let $ \mu_\bigdot $ be as above.  Then  $ \Phi \big(\overline{A(\mu_\bigdot)}\big) = \conv(\mu_\bigdot) $, or $ A(\mu_\bigdot) = \emptyset $.
\end{Lemma}
\begin{proof}
Let $ X = \overline{A(\mu_\bigdot)}$.  Assume that $ X $ is non-empty.   Let $P $ denote the moment polytope of $ X $.  We know that $ P $ is the convex hull of the set $ \{ \mu \in \cwl : t^\mu \in X \}$.    

For each $ w \in W $ consider the one parameter subgroup $ w\cdot \rho^\vee : \C^\times \rightarrow T $.  Let $ L \in A(\mu_\bigdot) $.  Since $ X $ is closed and $ T$-invariant, $ \lim_{s \rightarrow \infty} L \cdot (w \cdot \rho^\vee) (s) \in X $.  But since $ L \in S_w^{\mu_w} $, we see that $ \lim_{s \rightarrow \infty} L \cdot (w \cdot \rho^\vee) (s)  = t^{\mu_w} $.

Hence $ t^{\mu_w} \in X $ for all $ w \in W$.  Hence $ \conv( \mu_\bigdot ) \subset P $.

Conversely, if $ t^\nu \in X $, then $ t^\nu \in \overline{S_w^{\mu_w}} $ for each $ w \in W$.  So $ \nu \in C_w^{\mu_w} $. Hence $ \nu \in \cap_w C_w^{\mu_w} $.  Since $ \cap_w C_w^{\mu_w} = \conv(\mu_\bigdot) $ is convex, this shows that $ P \subset \conv(\mu_\bigdot) $.
\end{proof}

For each $ L \in \Gr $, let $ P(L) $ denote the pseudo-Weyl polytope corresponding to the GGMS stratum containing $ L $. 

\subsection{The functions $\D_\gamma $} \label{se:Dgam}
We now introduce constructible functions on the affine Grassmannian whose joint level sets are the GGMS strata.  These functions are new, but were motivated by the work of Speyer \cite{speyer}.  

If $ U $ is a vector space over $ \C $, the vector space $ U \otimes \K $ has a filtration by $ \cdots \subset U \otimes t \mathcal{O} \subset U \otimes \mathcal{O} \subset U \otimes t^{-1} \mathcal{O} \subset \cdots $.  We use this filtration to define a function $ \val $ on $ U \otimes \K $, by $ \val(u) = k $ if $ u \in U \otimes t^k \mathcal{O} $ and $ u \notin U \otimes t^{k+1} \mathcal{O} $.

Fix a high weight vector $ v_{\fund_i} $ in each fundamental representation $ V_{\fund_i} $ of $ G $.  For each chamber weight $ \gamma = w \cdot \fund_i $, let $ v_\gamma = \overline{w} \cdot v_{\fund_i} $. Since $ G $ acts on $ V_{\fund_i} $,  $ G(\K) $ acts on $ V_{\fund_i} \otimes \K $.

For each $ \gamma \in \Gamma $ define the function $ \D_\gamma $ by:
\begin{equation} \label{eq:Vdef}
\begin{aligned} 
\D_\gamma: \Gr & \rightarrow \Z \\
[g] &\mapsto \val ( g \cdot v_\gamma)
\end{aligned}
\end{equation}

This gives a well-defined function on $ \Gr = G(\mathcal{O}) \setminus G(\K) $, since if $ g \in G(\mathcal{O} ) $ and $ u \in V_{\fund_i} \otimes \K $, then $ \val( g \cdot u ) = \val( u) $.

The functions $ \D_\gamma $ have a simple structure with respect to the semi-infinite cells.  If $ \gamma = w \cdot \fund_i $, then $ v_\gamma $ is invariant under $ N_w(\K)$.  This immediately implies the following lemma.
\begin{Lemma} \label{th:SmDg}
Let $ w \in W $.
The function $ \D_{w \cdot \fund_i} $ takes the constant value $ \langle \mu, w \cdot \fund_i \rangle $ on $ S_w^\mu $.  In fact,
\begin{equation*}
S_w^\mu = \{ L \in \Gr : \D_{w \cdot \fund_i} (L) = \langle \mu, w \cdot \fund_i \rangle \text{ for all } i \}.
\end{equation*}
\end{Lemma}

Let $ M_\bigdot $ be a collection of integers, one for each chamber weight.  Then we consider the joint level set of the functions $ \D_\bigdot$,
\begin{equation}
A(M_\bigdot) := \{ L \in Gr : \D_\gamma(L) = M_\gamma \text{ for all } \gamma  \in \Gamma \}.
\end{equation}

Let $ \mu_\bigdot $ be a collection of coweights describing a pseudo-Weyl polytope.  Let $ M_\bigdot $ be the corresponding collection of integers defined by (\ref{eq:BZtoGGMS}).  Then by Lemma \ref{th:SmDg}, we have two descriptions of the GGMS stratum: $ A(\mu_\bigdot) = A(M_\bigdot)$. 

By Proposition \ref{th:BZtoGGMS}, we also have two different descriptions of the pseudo-Weyl polytope: $ \conv(\mu_\bigdot) = P(M_\bigdot) $. 

If the GGMS stratum is non-empty, then the GGMS stratum and the pseudo-Weyl polytope are related in two different ways:
\begin{gather*}
A(\mu_\bigdot) = A(M_\bigdot) = \{ L \in \Gr : P(L) = \conv(\mu_\bigdot) = P(M_\bigdot) \}, \\
\Phi \big(\overline{A(\mu_\bigdot)}\big) = \Phi \big(\overline{A(M_\bigdot)}\big) = \conv(\mu_\bigdot) = P(M_\bigdot),
\end{gather*}
where the first line of equations is by the definition of $P(L) $ and the second is by Lemma \ref{th:momentpoly}.

\section{BZ data and MV cycles} \label{se:BZdat}
Now we will give necessary and sufficient conditions on a collection $ M_\bigdot $ for $ \overline{A(M_\bigdot)} $ to be an MV cycle.  

\subsection{Generalized minors}
For this purpose, it is necessary to understand better the functions $ \D_\bigdot $.  To that end, we consider the \textbf{generalized minors} of Berenstein-Zelevinsky \cite{BZschub}.  For each chamber weight $ \gamma $ of level $ i $, they introduced the function
\begin{equation} \label{eq:mindef}
\begin{aligned}
\Delta_\gamma : G &\rightarrow \C \\
g &\mapsto \langle g \cdot v_\gamma , v_{-\fund_i} \rangle
\end{aligned}
\end{equation}
(note that $ v_{-\fund_i} \in V_{-w_0 \cdot \fund_i} = V_{\fund_i}^\dual $).

When $ G = SL_n $, a chamber weight of level $ i $ is just an $ i $ element subset of $ \{1, \dots, n\} $ and $ \Delta_\gamma (g) $ is the minor of $ g $ using the first $ i $ rows and column set $ \gamma $.

The function $ \D_\gamma $ on the affine Grassmannian is closely related to the valuation of $ \Delta_\gamma $.  In general, one can see that $ \val(\Delta_\gamma(g)) \ge \D_\gamma([g]) $ (see the remarks at the beginning of section \ref{se:off}).  We will show (in the course of the proof of Theorem \ref{th:BtoA}) that if $ L \in \Gr $, then there exists $ g \in G(\K) $ such that $ [g] = L $ and $ \D_\gamma(L) = \val(\Delta_\gamma(g)) $ for all $\gamma $.

\subsection{Tropical Pl\"ucker relations} \label{se:tpr}
Berenstein-Zelevinsky \cite{BZschub} established certain three-term Pl\"ucker relations among these generalized minors.  As our functions $ \D_\gamma $ are closely related to the valuation of these generalized minors, we would expect some relations among them coming from the tropical Pl\"ucker relations.  

In general, the process of passing from relations among Laurent series to relations among integers using $ \val $ is called tropicalization (see \cite{SS}).  Note that if $ a, b \in \K $, then
\begin{equation*}
\val(ab) = \val(a) + \val(b), \quad \val(a + b) \ge \min(\val(a) , \val(b)),
\end{equation*}
with equality holding in the second equation as long as the leading terms of $a $ and $ b $ do not cancel.  So if $ a,b,c,d \in \K $ satisfy the equation $ a = (b+c) d $, then the naive form of the tropicalization is
\begin{equation*}
A = \min(B,C) + D 
\end{equation*}
where $ A, B, C, D $ denote the valuations of $ a,b,c,d$.

We will show that this naive tropicalization is enough to understand the values of the $ \D_\gamma $ on an open subset of each MV cycle.  This motivates the following definition which originally appeared (though with a different motivation) in \cite{BZtpm}.

Let $ w \in W, i,j \in I $ be such that $ w s_i > w, ws_j > w $, $ i \ne j$.  We say that a collection $ \big( M_\gamma \big)_{\gamma \in \Gamma} $ satisfies the \textbf{tropical Pl\"ucker relation} at $ (w, i,j) $ if $a_{ij} = 0 $ or if
\begin{enumerate}
\item if $ a_{ij} = a_{ji} = -1 $, then
\begin{equation} \label{eq:A2trop}
M_{ws_i \cdot \fund_i} + M_{w s_j \cdot \fund_j} = \min(M_{w \cdot \fund_i} + M_{w s_i s_j \cdot \fund_j} , M_{w s_j s_i \cdot \fund_i} + M_{w \cdot \fund_j} ) ;
\end{equation}

\item if  $ a_{ij} = -1, a_{ji} = -2 $, then   
\begin{equation} \label{eq:B2trop1}
\begin{aligned}
M_{ws_j \cdot \fund_j} + M_{w s_i s_j \cdot \fund_j} + M_{ws_i \cdot \fund_i} = 
\min \Big( &2M_{w s_i s_j \cdot \fund_j} + M_{w \cdot \fund_i}, \\
&2M_{w \cdot \fund_j} + M_{w s_i s_j s_i \cdot \fund_i}, \\
&M_{w \cdot \fund_j} + M_{w s_j s_i s_j \cdot \fund_j} + M_{w s_i \cdot \fund_i} \Big), \\ 
M_{w s_j s_i \cdot \fund_i} + 2M_{w s_is_j \cdot \fund_j} + M_{w s_i \cdot \fund_i} = 
\min \Big(
&2M_{w \cdot \fund_j} + 2M_{w s_i s_j s_i \cdot \fund_i}, \\
&2M_{w s_j s_i s_j \cdot \fund_j} + 2M_{w s_i \cdot \fund_i}, \\
 &M_{w s_i s_j s_i \cdot \fund_i} + 2 M_{w s_i s_j \cdot \fund_j} + M_{w \cdot \fund_i} \Big)  ;
\end{aligned}
\end{equation}

\item if $ a_{ij} = -2, a_{ji} = -1 $, then 
\begin{equation} \label{eq:B2trop2}
\begin{aligned}
M_{w s_j s_i \cdot \fund_i} + M_{w s_i \cdot \fund_i} + M_{w s_i s_j \cdot \fund_j} = 
\min \Big( &2M_{w s_i \cdot \fund_i} + M_{w s_j s_i s_j \cdot \fund_j}, \\
&2M_{w s_i s_j s_i \cdot \fund_i} + M_{w \cdot \fund_j}, \\
&M_{w s_i s_j s_i \cdot \fund_i} + M_{w \cdot \fund_i} + M_{w s_i s_j \cdot \fund_j} \Big), \\
M_{w s_j \cdot \fund_j} + 2 M_{w s_i \cdot \fund_i} + M_{w s_i s_j \cdot \fund_j} = 
\min \Big( &2M_{w s_i s_j s_i \cdot \fund_i} + 2M_{w \cdot \fund_j}, \\
&2M_{w \cdot \fund_i} + 2M_{ w s_i s_j \cdot \fund_j}, \\
&M_{w \cdot \fund_j} + 2M_{w s_i \cdot \fund_i} + M_{w s_j s_i s_j \cdot \fund_j} \Big).
\end{aligned}  
\end{equation}
\end{enumerate}

We say that a collection $ M_\bigdot = \big( M_\gamma \big)_{\gamma \in \Gamma} $ satisfies the \textbf{tropical Pl\"ucker relations} if it satisfies the tropical Pl\"ucker relation at each $ (w, i, j) $. 

\subsection{BZ data}
A collection $ \big( M_\gamma \big)_{ \gamma\in \row } $ is called a \textbf{BZ (Berenstein-Zelevinsky) datum} of coweight $(\mu_1, \mu_2) $ if:
\begin{enumerate}
\item $ M_\bigdot $ satisfies the tropical Pl\"ucker relations.
\item $ M_\bigdot $ satisfies the edge inequalities (\ref{eq:nondeg}).
\item $ M_{\fund_i} = \langle \mu_1, \fund_i \rangle $ and $ M_{w_0 \cdot \fund_i} = \langle \mu_2, w_0 \cdot \fund_i \rangle $ for all $i$.
\end{enumerate}

The corresponding pseudo-Weyl polytope $ P(M_\bigdot) $ will have lowest vertex $ \mu_e = \mu_1 $ and highest vertex $ \mu_{w_0} = \mu_2$.  

Our main result, which will be proven in sections \ref{se:decomp} and \ref{se:overlap}, is the following characterization of MV cycles and polytopes.

\begin{Theorem} \label{th:BZcycle}
Let $ M_\bigdot $ be a BZ datum of coweight $ (\mu_1, \mu_2) $.
Then $\overline{A(M_\bigdot)} $ is an MV cycle of coweight $ (\mu_1, \mu_2) $, and each MV cycle arises this way for a unique BZ datum $ M_\bigdot $.

Hence a pseudo-Weyl polytope $ P(M_\bigdot) $ is an MV polytope if and only if $ M_\bigdot $ satisfies the tropical Pl\"ucker relations.
\end{Theorem}

In general if $ Y \subset  X $ is irreducible and $ f : X \rightarrow S $ is a constructible function, then there is a unique value $ s \in S $ such that $f^{-1}(s) \cap Y $ is dense in $ Y $.  In this situation, $ s $ is called the \textbf{generic value} of $ f $ on $ Y $.  

Using this language, Theorem \ref{th:BZcycle} says that if $ A $ is an MV cycle and if $ M_\gamma $ is the generic value of $ \D_\gamma $ on $A$ for each $ \gamma$, then $ M_\bigdot $ is a BZ datum.

\subsection{MV polytopes} \label{se:MVpoly}
In the case of $ G = SL_3 $, it is possible to give a very explicit description of the BZ data and MV polytopes.  In this case we have $ \Gamma = \{ 1, 2, 3, 12, 13, 23 \} $ where we use $ 2 $ as shorthand for $ (0,1,0) \in \wl $, $23$ for $(0,1,1) $, etc.

There is only one tropical Pl\"ucker relation (which occurs at $ (w= 1, i=1, j=2) $),
\begin{equation} \label{eq:plusl3}
M_2 + M_{13} = \min \{ M_1 + M_{23}, M_3 + M_{12} \}.
\end{equation}

Translated into the world of polytopes, we note that pseudo-Weyl polytopes for $ SL_3 $ are hexagons with every pair of opposite sides parallel and all sides meeting at $ 120^\circ $.  The above relation (\ref{eq:plusl3}) shows that a pseudo-Weyl  polytope is an MV polytope if and only if the distance between the middle pair of opposite sides is the maximum of the distances between the other two pairs of opposite sides.  Hence there are two possible forms for $ SL_3 $ MV polytopes, depending on which distance achieves this maximum.  Here are example of each of the two kinds (where $ \mu_1 $ marks the $ e$ vertex and $ \mu_2 $ marks the $ w_0 $ vertex).

\begin{figure}[h]
\begin{center}
\psfrag{0}{$\mu_1$} \psfrag{m}{$\mu_2$}
\epsfig{file=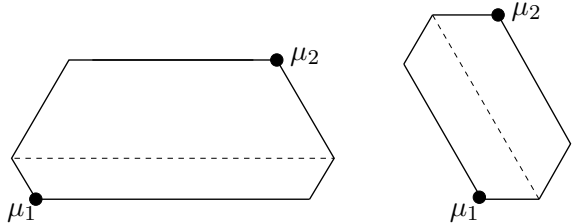,height=3cm}
\caption{The $SL_3 $ MV polytopes.}
\label{fig:sl3}
\end{center}
\end{figure}

In the case of $ G = Sp_4 $, there are two equivalent tropical Pl\"ucker relation (occurring at $ (w= 1, i = 1,j = 2) $ and at $ (w =1, i =2, j=1) $).  Examining the possible cases in either (\ref{eq:B2trop1}) or (\ref{eq:B2trop2}) shows there are the following four possible types of polytopes.

\begin{figure}[h]
\begin{center}
\psfrag{0}{$\mu_1$} \psfrag{m}{$\mu_2$}
\epsfig{file=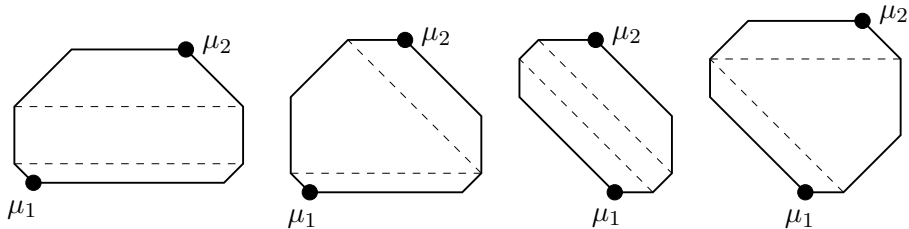,height=2.9cm}
\caption{The $Sp_4$ MV polytopes.}
\label{fig:sp4}
\end{center}
\end{figure}

In each case, there are certain interior edges going in root directions and connecting pairs of vertices.  These edges are shown by dotted lines in the above pictures.  The vertices connected by these edges are never the highest or lowest vertices.  Also, the edges in a particular picture do not cross.  Moreover, we see that for $ SL_3 $ and $ Sp_4 $ there is a 1-1 correspondence between maximal such arrangements and types of MV polytopes.  We do not have a good theoretical explanation for this phenomenon.   

Each tropical Pl\"ucker relation concerns the placement of the hyperplanes incident to a particular 2-face of the pseudo-Weyl polytope.  Hence we see that if $ \operatorname{rank}(G) > 2 $, then a pseudo-Weyl polytope is an MV polytope if and only if all of its 2-faces are MV polytopes.  So a pseudo-Weyl polytope is an MV polytope if and only if all of its 2-faces are rectangles (the MV polytopes for $ SL_2 \times SL_2 $) or one of the above types.  More generally, this shows that any face of an MV polytope is an MV polytope.  It is possible to understand this fact using the restriction map $ q_J $ introduced by Braverman-Gaitsgory \cite{BG} and further discussed \cite{MVcrystal}.

A small caveat is in order.  Each MV polytope comes with a labelling of its vertices by Weyl group elements.  When we look at a face of an MV polytope, this induces a labelling of the elements of that face by the corresponding Weyl group.  On the other hand, this labelling is automatic, because it is the only labelling consistent with its presentation as a pseudo-Weyl polytope (for example the ``$e$'' vertex always has to be the lowest weight vertex).  When we say that a face is an MV polytope, we really mean that we are considering this face along with its induced labelling.  This is important because as can be seen from $ SL_3 $ MV polytopes, the rotation/reflection of an MV polytope is not necessarily an MV polytope.

\section{Lusztig data decomposition} \label{se:decomp}

\subsection{Reduced words and paths}
Fix a reduced word $ \textbf{i} = (i_1, \dots, i_p) $ for an element $ w \in W $.   The word $ \textbf{i} $ determines a sequence of distinct Weyl group elements $ w^\textbf{i}_k := s_{i_1} \cdots s_{i_k} $ and distinct positive coroots $ \beta_k^\mathbf{i} := w^\textbf{i}_{k-1} \cdot \alpha^\vee_{i_k} $, $ k = 1 \dots p $.  In particular, when $ w = w_0 $, we get all the positive coroots this way. 

We say that a chamber weight $ \gamma $ is an $ \textbf{i} $-\textbf{chamber weight} if it is of the form $ w_k^\textbf{i} \cdot \fund_j $ for some $k, j$.  We write $ \Gamma^\textbf{i} $ for the set of all $ \textbf{i} $-chamber weights.  Let $ \gamma_k^\textbf{i} = w^\textbf{i}_k \cdot \fund_{i_k} $.  

Because of the relationship $ s_j \cdot \fund_i = \fund_i $ for $ j \ne i $, it is fairly easy to see that $ \Gamma^\wi $ consists of $ m + r $ elements: the $ \gamma_k^\wi $ and the fundamental weights (see \cite[Prop 2.9]{BZschub}).

It is worth keeping in mind the polytope combinatorics associated to this choice of reduced word.  Let $ \Sigma := W_{\rho^\vee} $ be the $ \rho^\vee$-Weyl polytope.  We will refer to this polytope as the \textbf{permutahedron}.  For each $ w \in W $, it has a vertex $ w w_0 \cdot \rho^\vee $ which we call the $ w$ vertex of $ \Sigma $.  For each $ w \in W $ and $ i \in I $, there is an edge connecting the $ w$ vertex and the $ws_i $ vertex.  Understanding the faces of the permutahedron is enough to understand the faces of any pseudo-Weyl polytope since there is a map from the set of faces of the permutahedron onto the set of faces of any pseudo-Weyl polytope (see appendix \ref{se:appendix}).  

A reduced word $ \wi $ determines a distinguished path $ w_0^\wi=e, w_1^\wi= s_{i_1}, w_2^\wi, \dots, w_m^\wi=w $ through the 1-skeleton of $ \Sigma $.  The $ k$th leg of this path is the vector $ w_{k-1}^\wi \cdot \rho - w_{k}^\wi \cdot \rho = \beta^\wi_k $.  The $ \wi$-chamber weights are exactly those dual to hyperplanes incident to the vertices along this path.    

\begin{Example} \label{eg:2}
Consider $ G = SL_3 $.  Let $ \wi = (1,2,1) $, then 
\begin{equation*}
w^\wi_1 = 213, \ w^\wi_2 = 231, \ w^\wi_3 = 321,
\end{equation*}
and
\begin{equation*}
\beta^\wi_1 = (1,-1,0), \ \beta^\wi_2 = (1,0,-1), \ \beta^\wi_3 = (0,1,-1).
\end{equation*}

Also,
\begin{equation*}
\gamma_1^\wi = 2, \ \gamma_2^\wi = 23, \ \gamma_3^\wi = 3,  
\end{equation*}
where we write $ (0,1,0) $ as $ 2$, $ (0,1,1) $ as $ 23 $, etc.

The fundamental weights $ 1, 12 $ are also $\wi$-chamber weights, so in fact every chamber weight is a $ \wi$-chamber weight except for $ 13$.

In figure \ref{fig:permut}, we show the permutahedron for $ SL_3$ along with the distinguished path corresponding to $ \wi $ and the hyperplanes defined by each chamber weight.

\begin{figure}  
\begin{center}
\psfrag{b1}{$\beta_1$}\psfrag{b2}{$\beta_2$}\psfrag{b3}{$\beta_3$}
\epsfig{file=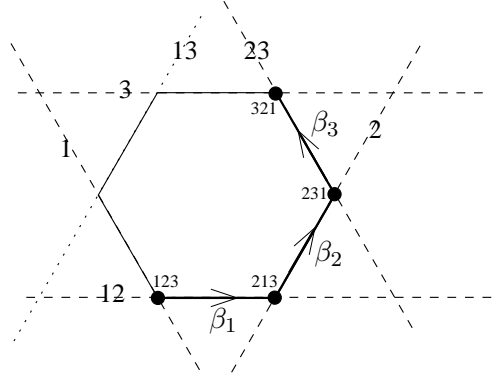,height=5cm}
\caption{The permutahedron for $ SL_3 $.}
\label{fig:permut}
\end{center}
\end{figure}
\end{Example}

\subsection{Lusztig data}
Let $ \wi = (i_1, \dots, i_m) $ be a reduced word for $ w_0 $.  If $ P = \conv(\mu_\bigdot) $ is a pseudo-Weyl polytope, we also get a distinguished path $ \mu_e, \mu_{s_{i_1}}, \mu_{s_{i_1} s_{i_2}}, \dots, \mu_{w_0} $ through the 1-skeleton of $ P $.  Let $ n_1, \dots, n_m $ be the sequence of lengths of the edges of this path.  We call the vector $ (n_1, \dots, n_m ) $ the $\wi$-\textbf{Lusztig datum} of P.  

Let $ n_\bigdot \in \N^m $.  We say that $ n_\bigdot $ is an $\wi$-\textbf{Lusztig datum} of coweight $ \mu $ if $ \mu = \sum_k n_k \beta^\wi_k $.  For such $ n_\bigdot$, let $ \mathcal{Q}^\wi(n_\bigdot) $ be the collection of pseudo-Weyl polytopes $ P = \conv(\mu_\bigdot) $, such that P has $ \wi$-Lusztig datum $ n_\bigdot $ and lowest vertex $ \mu_e = 0 $.  Note that if $ P \in \mathcal{Q}^\wi(n_\bigdot) $, then $ \mu_{w_0} = \sum_k n_k \beta^\wi_k = \mu $ is the coweight of the $\wi$-Lusztig datum of $ P $.

\begin{Example} \label{eg:2.5}
Continuing as in Example \ref{eg:2}, we see that there are three pseudo-Weyl polytopes with $ \wi $-Lusztig datum (2,1,1) and lowest vertex 0.  We will show that only the middle one is an MV polytope.\samepage
\begin{figure}[h]
\begin{center}
\epsfig{file=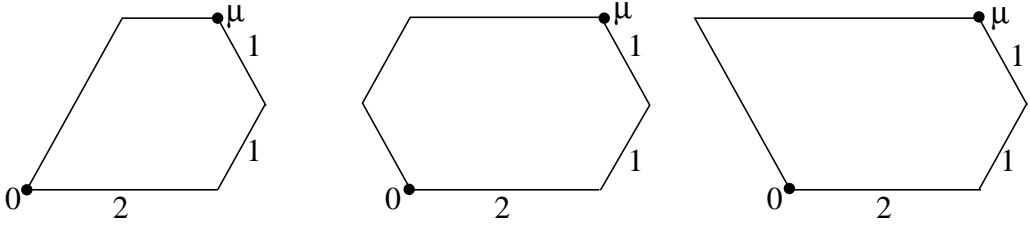, height=3cm}
\end{center}
\caption{The pseudo-Weyl polytopes with $ \wi$-Lusztig datum $(2,1,1) $.}
\end{figure}
\end{Example}

\subsection{The decomposition} 
With these considerations in mind, we proceed to discuss the decomposition according to Lusztig data.  Fix a reduced word $ \textbf{i} $ for $ w_0 $ and a coweight $ \mu \ge 0 $.  Let
\begin{equation*}
X(\mu) := S_e^0 \cap S_{w_0}^\mu.
\end{equation*}

Let $ A^\wi(n_\bigdot) := \{ L \in X(\mu) : P(L) \in \mathcal{Q}^\wi(n_\bigdot) \} $.  Since each pseudo-Weyl polytope has some $\wi$-Lusztig datum, we immediately have the following decomposition of $ X(\mu) $ into locally closed subsets.

\begin{Proposition} \label{th:Lddec}
\begin{equation*}
X(\mu) = \bigsqcup A^\wi(n_\bigdot) 
\end{equation*}
where the union is over all $\wi$-Lusztig data $n_\bigdot $ of coweight $ \mu $.
\end{Proposition}

Fix an $ \wi$-Lusztig datum $n_\bigdot $ of coweight $ \mu $.  Let $ \mu_k = \sum_{l=1}^k n_l \beta_l^\wi $.  Suppose that $ P $ is a pseudo-Weyl polytope with $\wi$-Lusztig datum $n_\bigdot $.  Then the $ w_k^\wi $ vertex of $ P $ is at position $ \mu_k $.  So if $ L \in A^\wi(n_\bigdot) $, then $ L $ lies in a GGMS stratum $ A(\nu_\bigdot) $ with $ \nu_{w_k^\wi} = \mu_k $.  This shows that
\begin{equation*}
A^\wi(n_\bigdot) = \bigcap_k S_{w^\wi_k}^{\mu_k}.
\end{equation*}

Let $ M_{\gamma_k^\wi} = \langle \mu_k, \gamma_k^\wi \rangle$.  Then by the length formula (\ref{eq:length}), we see that  $ \big( M_\gamma \big)_{\gamma \in \Gamma^\wi} $ is the unique solution to the system of equations
\begin{equation} \label{eq:Mfromn}
\begin{gathered}
n_k = - M_{w_{k-1}^\wi \cdot \fund_{i_k}} - M_{w_k^\wi \cdot \fund_{i_k}} - \sum_{j \ne i_k} a_{j,i_k} M_{w_k^\wi \cdot \fund_j} \text{ for all } k, \\
M_{\fund_i} = 0 \text{ for all } i.
\end{gathered}
\end{equation}

This system is upper triangular (note that each $ M_{\gamma_k^\wi} $ shows up for the first time in the equation with $ n_k $ on the left hand side) and so such a solution exists and is unique.  The solution is given by
\begin{equation} \label{eq:Mfromninverse}
M_{\gamma_k^\wi} = \sum_{l \le k} \langle \beta_l^\wi, \gamma_k^\wi \rangle n_l.
\end{equation}
For the proof, see Theorem 4.3 in \cite{BZschub}.

By Lemma \ref{th:SmDg} it follows that
\begin{equation} \label{eq:Andef}
A^\wi(n_\bigdot) = \{ L \in \Gr : \D_\gamma(L) = M_\gamma \text{ for all $ \wi $-chamber weights } \gamma \}.
\end{equation}

\begin{Example} \label{eg:3}
Continuing as in Example \ref{eg:2}, we see that in this case
\begin{gather*}
\mu_1 = (n_1, -n_1, 0), \ \mu_2 = (n_1 + n_2, -n_1, -n_2), \ \mu_3 = (n_1 + n_2, n_3 - n_1, -n_2 - n_3), \\
n_1 = -M_2, \ n_2 = -M_{23} + M_2, \ n_3 = -M_2 - M_3 + M_{23}.
\end{gather*}
\end{Example}

The goal of this section is to prove the following  result.

\begin{Theorem} \label{th:Ldirr}
For each $ \wi$-Lusztig data of coweight $ \mu$, $ \overline{A^\wi(n_\bigdot)} $ is an irreducible component of $ \overline{X(\mu)} $.  Moreover each component of $ \overline{X(\mu)} $ appears exactly once this way.
\end{Theorem}

The following elementary algebraic geometry lemma will prove quite useful.

\begin{Lemma} \label{th:alggeo}
Let $ X $ be a reducible algebraic set with $ n $ components.  Suppose that $ X = \sqcup C_k $ is a decomposition of $ X $ into $ n $ irreducible constructible subvarieties. Then $ \overline{C_1}, \dots, \overline{C_n} $ are the distinct irreducible components of $ \overline{X} $.
\end{Lemma}

\begin{proof}
Let $ A_1, \dots, A_n $ denote the irreducible components of $ \overline{X} $.  Then $ \overline{X} = \cup \overline{C_i} $ , so
\begin{equation*}
A_j = \bigcup A_j \cap \overline{C_i}.
\end{equation*}  
Since $ A_j $ is irreducible and each $ A_j \cap \overline{C_i} $ is closed, $ A_j = A_j \cap \overline{C_i} $ for some $ i$.  So $ A_j \subset \overline{C_i} $.  By similar reasoning, there exists $ k $ such that $ \overline{C_i} \subset A_k $.  Hence $ A_j \subset \overline{C_i} \subset A_k $.  Since the $ A_j $ are the components, each listed once, $ j = k $ and so $ A_j = \overline{C_i} $.  Continuing this argument shows that there exists a map $ \sigma $ of $ \{ 1, \dots, n \} $ to $ \{ 1, \dots, n \} $ such that $ A_j = \overline{C_{\sigma(j)}} $.  This map is injective since the $ A_j $ are distinct.  Hence it is bijective as desired. 
\end{proof}
 
The number of $ \wi $-Lusztig data of coweight $ \mu $ is $\kpf(\mu) $ which equals the number of components of $ \overline{X(\mu)} $.  So to prove Theorem \ref{th:Ldirr}, it suffices to show that $ A^\wi(n_\bigdot) $ is irreducible for each Lusztig datum $ n_\bigdot $.  To prove this, we will use another basic algebraic geometry fact, that the image of an irreducible variety is irreducible.  Hence our goal is to construct a surjective map from an irreducible variety onto $ A^\wi(n_\bigdot) $.  To that end, we will examine certain parametrizations of $ N $ introduced by Lusztig and Berenstein-Fomin-Zelevinsky.

\subsection{Parametrizations of $ N $}
Fix $ w \in W $.  Following Berenstein-Zelevinsky \cite{BZschub}, we will define the \textbf{twist} automorphism $ \eta_w : N \cap B_- w B_- \rightarrow N \cap B_- w B_- $.  First, let $ x \mapsto x^T $ be the involutive Lie algebra anti-automorphism of $ \mathfrak{g} $ given by
\begin{equation*}
E_i^T = F_i, \quad F_i^T = E_i, \quad H_i^T = H_i,
\end{equation*}
where $E_i, F_i, H_i $ denote the Chevalley generators of $ \mathfrak{g} $, corresponding to the maps $ \psi_i $ of $ SL_2 $ into $ G $.  We use the same notation $ g \mapsto g^T $ for the corresponding involutive anti-automorphism of $ G $.  

For $ y \in N \cap B_- w B_-$, we define $ \eta_w(y) $ to be the unique element in the intersection $ N \cap B_- \overline{w} y^T $. See \cite{BZschub} for proof that this function is well-defined.

We define $ \mathbf{x}_i : \C \rightarrow N $ by
\begin{equation*}
\mathbf{x}_i(a) = \psi_i \Big(
\begin{bmatrix}
1&  a \\
0 & 1 \\
\end{bmatrix} \Big).
\end{equation*}

Let $ \wi $ be a reduced word for $ w $ and let $ p $ be its length.  Following \cite{Lbookbook, BZschub}, we define regular maps $ \mathbf{x}_\wi $ and $ \mathbf{y}_\wi $ from $ (\C^\times)^p $ to $ N $,
\begin{align*}
\mathbf{x}_\wi(b_1, \dots, b_p) &= \mathbf{x}_{i_p}(b_p) \cdots \mathbf{x}_{i_1}(b_1), \\
\mathbf{y}_\wi(b_1, \dots, b_p) &= \eta_{w^{-1}}^{-1}(\mathbf{x}_\wi(b_1, \dots, b_p)).
\end{align*}

Berenstein-Fomin-Zelevinsky established the following result, which they call the \textbf{Chamber Ansatz}, which provides an inverse for $ y $.

\begin{Theorem} \label{th:ChAns}
Let $ y = \mathbf{y}_\wi(b_1, \dots, b_p) $. Then
\begin{equation} \label{eq:ChAns}
b_k = \frac{1}{\Delta_{w^\wi_{k-1} \cdot \fund_{i_k}}(y) \Delta_{w^\wi_k \cdot \fund_{i_k}}(y) } \prod_{j \ne i_k} \Delta_{w^\wi_k \cdot \fund_j}(y)^{-a_{j, i_k}} \text{ for all } k.
\end{equation}

Conversely, $\Delta_\gamma(y) $ is a monomial in the $ b_k $ whenever $ \gamma $ is a $\wi$-chamber weight.

Moreover, if $ w = w_0 $, then $ \mathbf{y}_\wi $ is a biregular isomorphism onto $ \{ y \in N : \Delta_\gamma(y) \ne 0 $ for all $\wi$-chamber weights $ \gamma  \} $.
\end{Theorem}

\begin{proof}
The first part of this theorem is exactly Theorem 1.4 in \cite{BZschub} and Theorem 2.19 in \cite{FZ}, except that we have reversed the order of the reduced word.

The system (\ref{eq:ChAns}) is the same as the system (\ref{eq:Mfromn}), except it is written multiplicatively instead of additively.  We have already observed that (\ref{eq:Mfromn}) is invertible, hence so is (\ref{eq:ChAns}) and so $ \Delta_\gamma(y) $ is a monomial in the $ b_k$ (this monomial is given in additive form in (\ref{eq:Mfromninverse})).

To prove the last statement, let $ U = \{ y \in N : \Delta_\gamma(y) \ne 0 $ for all $\wi$-chamber weights $ \gamma \} $.  The first half of the theorem provides a map $ U \rightarrow (\C^\times)^m $ which is a left inverse to $ \mathbf{y}_\wi $. Hence $ \mathbf{y}_\wi $ is injective.  

So it suffices to show that $ \mathbf{y}_\wi $ is surjective.  Let $ y \in U $ and determine $ b_k $ from $ y $ by (\ref{eq:ChAns}).  Let $ y' = \mathbf{y}_\wi(b_\bigdot) $.  By the above observations, the generalized minors $ \Delta_\gamma $ take the same values on $ y, y' $ for each $\wi$-chamber weight $ \gamma $.  But by the results of \cite{BZschub}, every function on $ N $ is a rational function of the $ \Delta_\gamma $ for $ \gamma $ an $\wi$-chamber weight.  Hence every function on $ N $ takes the same values on $ y' $ and $ y $.  Since $ N $ is affine, this shows that $ y' = y $ and so $ \mathbf{y}_\wi $ is surjective.
\end{proof}

\begin{Example}
We continue from Example \ref{eg:3}.  In this case:
\begin{equation*}
\mathbf{x}_\wi(b_1,b_2, b_3) = 
\begin{bmatrix}
1 & b_1 + b_3 & b_2b_3 \\
0 & 1 & b_2 \\
0 & 0 & 1 \\
\end{bmatrix} \ \text{ and } \  
\mathbf{y}_\wi(b_1, b_2, b_3) = 
\begin{bmatrix}
1 & \frac{1}{b_1} & \frac{1}{b_2 b_3} \\
0 & 1 & \frac{b_1 + b_3}{b_2 b_3} \\
0 & 0 & 1 \\
\end{bmatrix}.
\end{equation*}

So
\begin{equation*}
b_1 = \frac{1}{\Delta_2(y)}, \ b_2 = \frac{\Delta_2(y)}{\Delta_{23}(y)}, \ b_3 = \frac{\Delta_{23}(y)}{\Delta_2(y)\Delta_3(y)} 
\end{equation*}
as claimed in Theorem \ref{th:ChAns}.
\end{Example}

Note that the map $ \mathbf{y}_\wi $ is a map of varieties over $ \C $.  By the formal loop space functor, there is a corresponding map of ind-schemes over $ \C $,  $ \K^p  \rightarrow  N(\K) $.  Moreover, the obvious analogue of Theorem \ref{th:ChAns} holds in this setting.

\subsection{Mapping onto the MV cycles}

Fix a reduced word $ \wi $ for $ w_0 $.  Let $ n_\bigdot $ be a Lusztig datum of coweight $ \mu $. Let $ \big( M_\gamma \big)_{\gamma \in \Gamma^\wi} $ be determined from the $ n_\bigdot $ by (\ref{eq:Mfromn}).  

Let 
\begin{equation*}
B(n_\bigdot) := \{ (b_1, \dots, b_m) \in \K^m : \val(b_k) = n_k \text{ for all } k \}.
\end{equation*}

The goal of the rest of this section is to prove the following theorem.

\begin{Theorem} \label{th:BtoA}
If $ b_\bigdot \in B(n_\bigdot) $, then $ [\mathbf{y}_\wi(b_\bigdot)] \in A^\wi(n_\bigdot) $.  Moreover, each $ L \in A^\wi(n_\bigdot) $ has a representative of the form $ \mathbf{y}_\wi(b_\bigdot) $ for some $ b_\bigdot \in B(n_\bigdot) $.  Hence the restriction of $ \mathbf{y}_\wi $ to $ B(n_\bigdot) $ combined with the surjection $ G(\K) \rightarrow \Gr $ provides a surjective morphism $ B(n_\bigdot) \rightarrow A^\wi(n_\bigdot) $.
\end{Theorem}

Note that $ B(n_\bigdot) $ is irreducible, since it is isomorphic to a product of $ m $ copies of $ \C^\times $ and $ m $ copies of $ \mathcal{O} $.  Hence by the remarks following Lemma \ref{th:alggeo}, proving Theorem \ref{th:BtoA} will complete the proof of Theorem \ref{th:Ldirr}.

As a first step towards Theorem \ref{th:BtoA}, we establish the following lemma.

\begin{Lemma} \label{th:iff}
Let $ b_\bigdot \in \K^m$.
Let $ y = \mathbf{y}_\wi(b_\bigdot) $.

Then $ b_\bigdot \in B(n_\bigdot) $ if and only if $ \val(\Delta_\gamma(y)) = M_\gamma $ for all $\wi$-chamber weights $ \gamma $.
\end{Lemma}

\begin{proof}
By Theorem \ref{th:ChAns}, we see that
\begin{equation}
\val(b_k) = - \val(\Delta_{w^\wi_{k-1} \cdot \fund_{i_k}}(y)) - \val(\Delta_{w^\wi_k \cdot \fund_{i_k}}(y)) - \sum_{j \ne i_k} a_{j, i_k} \val(\Delta_{w_k^\wi \cdot \fund_j} (y))
\end{equation}
for all $ k $.  Also since $y \in N(\K) $, $ \Delta_{\fund_i}(y) = 1 $ and so $ \val(\Delta_{\fund_i}(y)) = 0 $ for all $ i $. 

This is the same system of equations as (\ref{eq:Mfromn}), with $\val(b_k) $ instead of $ n_k $ and $ \val(\Delta_\gamma(y)) $ instead of $ M_\gamma $.  Since (\ref{eq:Mfromn}) is an invertible linear system, this shows that $ \val(b_k) = n_k $ for all $ k $ if and only if $ \val(\Delta_\gamma(y)) = M_\gamma $ for all $ \wi$-chamber weights $ \gamma $.
\end{proof}

\subsection{Off-minors} \label{se:off}
To complete the proof of Theorem \ref{th:BtoA}, we will need a further investigation of the relation between $ \D_\gamma $ and the valuation of $ \Delta_\gamma$.

Let $ U $ be a finite-dimensional vector space over $ \C $.  Earlier, we defined a function $ \val : U \otimes \K \rightarrow \Z $.  Note that if $ u \in U \otimes \K $, then
\begin{equation*}
\val(u) = \min_{\xi \in U^\dual} \val (\langle u, \xi \rangle) ,
\end{equation*}
where on the right, $ \val $ denotes the usual valuation map on $ \K $. In fact, it is enough to take the min over a basis for $ U^\dual $.

Let us apply the above result to our situation.  We see that if $ \gamma $ is a chamber weight of level $ i $, then 
\begin{equation} \label{eq:valismin}
\D_\gamma([y]) = \val( y \cdot v_\gamma ) = \min_{\xi \in V_{\fund_i}^\dual} \val(\langle y \cdot v_\gamma , \xi \rangle). 
\end{equation}
In particular, $ \xi = v_{- \fund_i} $ shows up on the right hand side and so $ \val(\Delta_\gamma(y)) $ appears in the min (see (\ref{eq:mindef})).  We would like to show that the minimum is attained there when $ y = \mathbf{y}_\wi(b_\bigdot) $ and $ b_\bigdot \in B(n_\bigdot) $.

Using a Bruhat decomposition of $ G(\K) $ it is possible to show that we need to take only extremal weight vectors $ \xi $ in the $ \min $ above.  However, we will not need this. 

We call $ \langle y \cdot v_\gamma, \xi \rangle $  an \textbf{off-minor} of $ y $.  In the case $ G= SL_n $ and $ \xi = v_\delta $, it is the minor of $ y $ using $ \gamma $ as the set of columns and $ \delta $ as the set of rows (where we index the usual basis for $ V_{\fund_i}^\dual $ by the $ i $ element subsets of $\{ 1, \dots, n \} $).

The following lemma is a generalization of Lemma 3.1.3 from \cite{BFZA}, which dealt with the case $ G = SL_n $.
\begin{Lemma} \label{th:evaloff}
Let $ w \in W $.  Let $ \xi \in V_{\fund_i}^\dual $.  Let $x \in N \cap B_- w^{-1} B_- $ and $ y = \eta_{w^{-1}}(x) $.  Then
\begin{equation*}
\frac{\langle y \cdot v_{w \cdot \fund_i} , \xi \rangle}{\Delta_{w \cdot \fund_i} (y)} = \frac{\langle x^T \cdot v_{\fund_i}, \xi \rangle}{\langle v_{\fund_i}, v_{-\fund_i} \rangle}
\end{equation*}
\end{Lemma}

\begin{proof}
Since $ x = \eta_{w^{-1}}(y) $, there exists $ p \in N_- $ and $ d \in T $ such that $ pdx = \overline{w^{-1}} y^T $.  Note that $ \overline{w}^{-1} = \overline{w}^T $ (since $ \overline{s_i}^{-1} =\overline{s_i}^T $ by an $SL_2 $ calculation). 
Hence, $ y = x^T d^T p^T \overline{w^{-1}}$, and so  
\begin{equation} \label{eq:offmin1}
\langle y \cdot v_{w \cdot \fund_i}, \xi \rangle = \langle x^T d^T p^T \overline{w^{-1}} \cdot v_{w \cdot \fund_i}, \xi \rangle = \fund_i(rd) \langle x^T \cdot v_{\fund_i}, \xi \rangle,
\end{equation}
where $ r = \overline{w^{-1}} \overline{w} \in T $.

Similarly,
\begin{equation} \label{eq:offmin2}
\langle y \cdot v_{w \cdot \fund_i}, v_{-\fund_i} \rangle = \fund_i(rd) \langle x^T \cdot v_{\fund_i}, v_{-\fund_i} \rangle = \fund_i(rd) \langle v_{\fund_i}, (x^T)^{-1} \cdot v_{-\fund_i} \rangle = \fund_i(rd) \langle v_{\fund_i}, v_{-\fund_i} \rangle
\end{equation}
since $ x^T \in N_- $, so $ (x^T)^{-1} \in N_- $ and hence $ (x^T)^{-1} \cdot v_{-\fund_i} = v_{-\fund_i} $.  

Taking the ratio of (\ref{eq:offmin1}) and (\ref{eq:offmin2}) gives the desired result.
\end{proof}

This result allows us to express certain off-minors of $ y $ in terms of $ x $.  To express them all, we will also need the following lemma of Berenstein-Zelevinsky.

\begin{Lemma}[{\cite[Proposition 5.4]{BZschub}}]
Let $ \wi $ be a reduced word for $ w_0 $, let $ 1 \le k \le m $, let $ w = w^\wi_k $, and let $y = \mathbf{y}_\wi(b_1, \dots, b_m) $.  Then $ y $ admits a factorization $ y = y' y'' $ where $ y' = \mathbf{y}_{(i_1, \dots, i_k)}(b_1, \dots, b_k) $, and $ y'' \in w N w^{-1} $.
\end{Lemma}

These last two lemmas combine in the following result describing the off minors.

\begin{Proposition} \label{th:offpol}
Let $\wi $ be a reduced word for $ w_0 $, let $ \xi \in V_{\fund_i}^\dual $, and let $ \gamma $ be an $ \wi $-chamber weight of level $ i $.  Let $ y = \mathbf{y}_\wi(b_1, \dots, b_m) $. Then
\begin{equation*}
\frac{\langle y \cdot v_\gamma , \xi \rangle}{\Delta_\gamma(y)} 
\end{equation*}
is a polynomial in the $ b_k $.
\end{Proposition}

\begin{proof}
Since $ \gamma $ is an $\wi$-chamber weight, $ \gamma = w^\wi_k \cdot \fund_i $ for some $ k $.  Let $ w = w^\wi_k $.  By the previous lemma, $ y = y' y'' $, where $ y' = \mathbf{y}_{(i_1, \dots, i_k)}(b_1, \dots, b_k) $ and $y'' \in w N w^{-1} $.

Then $ y \cdot v_\gamma = y' y'' \cdot v_\gamma = y' \cdot v_\gamma $ since $ \gamma = w \cdot \fund_i $ and $ y'' \in w N w^{-1} $.

So
\begin{equation*}
\frac{\langle y \cdot v_\gamma , \xi \rangle}{\Delta_\gamma(y)} = \frac{ \langle y' \cdot v_\gamma, \xi \rangle}{\Delta_\gamma(y')} = \frac{\langle {x'}^T \cdot v_{\fund_i} , \xi \rangle}{ \langle v_{\fund_i}, v_{-\fund_i} \rangle },
\end{equation*}
where $ x' = \mathbf{x}_{(i_1, \dots, i_k)}(b_1, \dots b_k) $.  The first equality is by the above analysis and the second is by Lemma \ref{th:evaloff}.

Any regular function of $ {x'}^T $ is a polynomial in the $ b_k $ (since the extension of $ \mathbf{x}_{(i_1, \dots, i_k)} $ to $\C^k $ is regular) and so the result follows.
\end{proof} 

We are now ready to prove Theorem \ref{th:BtoA}.

\begin{proof}[Proof of Theorem \ref{th:BtoA}]
First, we will show that if $ b_\bigdot \in B(n_\bigdot) $, then $ [\mathbf{y}_\wi(b_\bigdot)] \in A^\wi(n_\bigdot) $. 

Fix $ b_\bigdot \in B(n_\bigdot) $ and let $ y = \mathbf{y}_\wi(b_\bigdot) $.

By (\ref{eq:Andef}), $ [y] \in A^\wi(n_\bigdot) $ if $ \D_\gamma([y]) = M_\gamma $ for all $ \wi$-chamber weights $ \gamma $.

By Lemma \ref{th:iff}, $ \val(\Delta_\gamma(y)) = M_\gamma $.  So to prove that $ [y] \in A^\wi(n_\bigdot) $, it suffices to show that $ \val(\Delta_\gamma(y)) = \D_\gamma([y]) $.

By (\ref{eq:valismin}), it suffices to show that $ \val( \langle y \cdot v_\gamma, \xi \rangle ) \ge \val(\Delta_\gamma(y)) $ for any $ \xi \in V_{\fund_i}^\dual $.  By Proposition \ref{th:offpol}, 
\begin{equation*}
\frac{\langle y \cdot v_\gamma, \xi \rangle}{\Delta_\gamma(y)} = P(b_1, \dots, b_m)
\end{equation*}
for some polynomial $ P $.  But $ \val (b_k) = n_k \ge 0 $ for all $ k $, so $ \val(P(b_1, \dots, b_m)) \ge 0 $.  Hence $ \val( \langle y \cdot v_\gamma, \xi \rangle) - \val(\Delta_\gamma(y)) \ge 0 $ as desired.

So we conclude that $ [y] \in A^\wi(n_\bigdot) $, as desired.

Next, we need to check that if $ L \in A^\wi(n_\bigdot) $, then $ L = [\mathbf{y}_\wi(b_\bigdot)] $ for some $ b_\bigdot \in B(n_\bigdot) $.  Suppose we know that there exists $ y \in N(\K) $ such that $ L = [y] $ and $\val(\Delta_\gamma(y)) = \D_\gamma(L) $ for all $ \gamma $.  Then if $ \gamma $ is a $\wi$-chamber weight, by (\ref{eq:Andef}) $\D_\gamma(L) = M_\gamma $, so $\val(\Delta_\gamma(y)) = M_\gamma $.  In particular, $ \Delta_\gamma(y) $ is non-zero for all $\wi$-chamber weights $ \gamma $.  Hence by Theorem \ref{th:ChAns}, there exist $ (b_1, \dots, b_m) \in (\K^\times)^m $ such that $ y = \mathbf{y}_\wi(b_1, \dots, b_m) $.  By Lemma \ref{th:iff}, we see that $ b_k \in B(n_\bigdot) $ as desired.  Hence this completes the proof of the theorem.

So now we will prove the existence of $ y $ as above.  Since $ A(n_\bigdot) \subset S_e^0 $, $ L $ has a representative $ g \in N(\K) $.  Let $ h \in N(\C) $.  So $[h^{-1}g] = [g] = L $.  We would like to find $ h $ such that $ \val(\Delta_\gamma(h^{-1}g)) = \D_\gamma(L) $ for all chamber weights $ \gamma $.

Let $ \gamma $ be a chamber weight of level $ i $ and let $ d = \D_\gamma([g]) $.  Let $ u_1, \dots, u_N $ be a basis for $ V_{\fund_i} $ with dual basis $ u_1^\dual, \dots, u_N^\dual $ for $ V_{\fund_i}^\dual $.

Then
\begin{equation*}
\Delta_\gamma(h^{-1}g) = \langle h^{-1}g \cdot v_\gamma, v_{-\fund_i} \rangle = \langle g \cdot v_\gamma, h \cdot v_{-\fund_i} \rangle.
\end{equation*}

Now $ h \cdot v_{-\fund_i} = \sum_s c_s u_s^\dual $ for some $ c_s \in \C $.  Hence
\begin{equation*}
\Delta_\gamma(h^{-1}g) = \sum_s c_s \langle g \cdot v_\gamma, u_s^\dual \rangle.
\end{equation*}

Let $ p_s $ be the coefficient of $ t^d $ in $ \langle g \cdot v_\gamma, u_s^\dual \rangle $.  Since $ d = \D_\gamma([g]) = \min_s \val (\langle g \cdot v_\gamma, u_s^\dual \rangle) $, we see that $ p_s $ is nonzero for some $ s $.  Extracting the coefficient of $ t^d $ from the above equation shows that $ \val( \Delta_\gamma(h^{-1} g)) = d $ if and only if $ \sum_s p_s c_s \ne 0 $.  

Let $ u = \sum_s p_s u_s $.  Thus, we see that 
\begin{equation*}
\val(\Delta_\gamma(h^{-1}g)) = \D_\gamma(L) \text{ if and only if } \langle u, h \cdot v_{-\fund_i} \rangle \ne 0 
\end{equation*}

Note that $ h \mapsto \langle u, h \cdot v_{-\fund_i} \rangle $ is a non-zero regular function on $ N $, since $u \ne 0 $ and since $ V_{-\fund_i} $ is generated by $ N $ acting on $ v_{-\fund_i}$.  Similarly, for each $ \gamma \in \Gamma $, there is a non-zero regular function $ f_\gamma $ such that  $ \val(\Delta_\gamma(h^{-1}g)) = \D_\gamma(L) $ if and only if $ f_\gamma(h) \ne 0 $.  Since $ N $ is irreducible, the product of these non-zero functions is non-zero and so we can find $ h $ such that $ \val(\Delta_\gamma(h^{-1}g)) = \D_\gamma(L) $ for all $ \gamma $, as desired.
\end{proof}

\section{Piecing together} \label{se:overlap}
By Theorem \ref{th:Ldirr}, if $ \wi $ is a reduced word for $ w_0 $ and $ n_\bigdot $ is an $ \wi$-Lusztig datum, then $ A^\wi(n_\bigdot) $ is an MV cycle and all MV cycles arise this way.  So for each $ \wi $ we get a bijection from $\N^m$ to the set of MV cycles.  We call the inverse of this bijection the $\wi$-\textbf{Lusztig datum} of the MV cycle.   

To complete the proof of Theorem \ref{th:BZcycle}, it will be necessary to understand how the $\wi$-Lusztig datum varies when we change the reduced word $\wi$.  To do this, we will consider the overlap in the different decompositions of $ X(\mu) $ by Lusztig data.  

In this section, a reduced word will always mean a reduced word for $ w_0 $.

\subsection{Local picture} \label{se:loc}
Two reduced words $ \wi, \wi' $ are said to be related by a $d$-\textbf{braid move} involving $ i, j$, starting at position $ k $, if 
\begin{align*}
\wi &= (\dots, i_k, i, j, i, \dots, i_{k+d+1}, \dots), \\
\wi' &= (\dots, i_k, j, i, j, \dots, i_{k+d+1}, \dots), 
\end{align*}
where $ d $ is the order of $ s_i s_j $.

\begin{figure}
\begin{center}
\psfrag{wi}{$\wi$} \psfrag{wi'}{$\wi'$} \psfrag{wk}{$w_k$} \psfrag{wkd}{$w_{k+d}$} 
\psfrag{i}{$\scriptstyle{i}$} \psfrag{j}{$\scriptstyle{j}$}
\epsfig{file=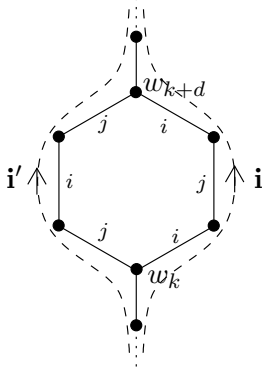,height=5cm}
\caption{Two reduced words related by a 3-braid move.}
\label{fig:twopaths}
\end{center}
\end{figure}

Recall that reduced words correspond to minimal length paths from the e vertex to the $w_0$ vertex of the permutahedron.  If $ \wi, \wi'$ are related as above, then $ w^\wi_l = w^{\wi'}_l $, for $ l \notin \{ k+1, \dots, k+d-1 \}$.  So the two paths agree for the first $ k $ vertices and then agree again at vertex $ k+d $ and later.  Moreover, the $w_l^\wi$ and $w_l^{\wi'}$ vertices for $ l \in \{ k, \dots, k+d \} $ all lie on the same 2-face of the permutahedron.  Namely, they lie on the 2-face which contains the $ w $ vertex and is dual to the chamber weights $ w \cdot \fund_p $ for $ p \ne i, j$, where $ w = w^\wi_k $.  This 2-face will be a 2d-gon (see Figure \ref{fig:twopaths}).  

Following Lusztig \cite{Lbookbook}, Berenstein-Zelevinsky studied the relationship between $ \mathbf{y}_\wi $ and $ \mathbf{y}_{\wi'} $.

\begin{Proposition}[{\cite[Theorem 3.1]{BZschub}}] \label{th:paramtrans}
Let $ \wi, \wi' $ be as above.
Suppose that $ \mathbf{y}_\wi(b_\bigdot) = \mathbf{y}_{\wi'}(b'_\bigdot) $.

For $ l \notin \{k+1, \dots, k +d\} $, $ b_l = b'_l $.  For other $ l $ we have the following case by case formulas. 
\begin{enumerate}
\item If $ a_{ij} = 0 $, then $ d =2 $ and 
\begin{equation*}
b'_{k+1} = b_{k+2}, \ b'_{k+2} = b_{k+1}.
\end{equation*}

\item  If $ a_{ij} = a_{ji} = -1$, then $ d = 3 $ and
\begin{equation} \label{eq:paramtrans2}
\begin{gathered}
b'_{k+1} = \frac{b_{k+2} b_{k+3} }{\pi}, \ b'_{k+2} = b_{k+1} + b_{k+3}, \ b'_{k+3} = \frac{b_{k+1} b_{k+2}}{\pi}, \\
\text{where } \ \pi = b_{k+1} + b_{k+3}.
\end{gathered}
\end{equation}

\item If $a_{ij} = -1, a_{ji} = -2$, then $ d = 4$ and
\begin{equation} \label{eq:paramtrans3}
\begin{gathered}
b'_{k+1} = \frac{b_{k+2} b_{k+3} b_{k+4} }{\pi_1}, \ b'_{k+2} = \frac{\pi_1^2}{\pi_2}, \ b'_{k+3} = \frac{\pi_2}{\pi_1}, \ b'_{k+4} = \frac{b_{k+1} b_{k+2}^2 b_{k+3}}{\pi_2}, \\
\text{where } \ \pi_1 = b_{k+1}b_{k+2} + (b_{k+1}+ b_{k+3})b_{k+4}, \ \pi_2 = b_{k+1}(b_{k+2} + b_{k+4})^2 + b_{k+3} b_{k+4}^2.
\end{gathered}
\end{equation}

\item If $ a_{ij} = -2, a_{ji} = -1$, then $ d = 4 $ and
\begin{equation} \label{eq:paramtrans4}
\begin{gathered}
b'_{k+1} = \frac{b_{k+2} b_{k+3}^2 b_{k+4} }{\pi_2}, \ b'_{k+2} = \frac{\pi_2}{\pi_1}, \ b'_{k+3} = \frac{\pi_1^2}{\pi_2}, \ b'_{k+4} = \frac{b_{k+1} b_{k+2} b_{k+3}}{\pi_1}, \\
\text{where } \ \pi_1 = b_{k+1}b_{k+2} + (b_{k+1}+ b_{k+3})b_{k+4}, \ \pi_2 = b_{k+1}^2b_{k+2} + (b_{k+1} + b_{k+3})^2 b_{k+4}.
\end{gathered}
\end{equation}

\end{enumerate}

Conversely, suppose that $ b_\bigdot \in (\C^\times)^m $ is such that the denominators in the above expressions do not vanish.  Define $ b'_\bigdot $ by the above expressions.  Then $\mathbf{y}_\wi(b_\bigdot) = \mathbf{y}_{\wi'}(b'_\bigdot) $.
\end{Proposition}

The first part of this proposition is directly from \cite{BZschub}.  The last statement follows from the same reasoning as in our proof of the second statement of Theorem \ref{th:ChAns}.  

Note that Proposition \ref{th:paramtrans} holds over $ \K $ as well.

\begin{Proposition} \label{th:valtrans}
Let $ n_\bigdot $ be an $\wi$-Lusztig datum of coweight $ \mu $.
Then there exists a non-empty open subset $ U $ of $ B(n_\bigdot) $ such that for each $ b_\bigdot \in U$, there exists $ b'_\bigdot \in \K^m $ such that $ \mathbf{y}_\wi(b_\bigdot) = \mathbf{y}_{\wi'}(b'_\bigdot) $  and the following formulas holds for $ n'_l := \val(b'_l)$.
\begin{enumerate}
\item If $a_{ij} = 0 $, then $ d =2 $ and 
\begin{equation*}
n'_{k+1} = n_{k+2}, \ n'_{k+2} = n_{k+1}.
\end{equation*}

\item If $ a_{ij} = a_{ji} = -1$, then $ d = 3 $ and
\begin{equation}  \label{eq:valtrans2}
\begin{gathered}
n'_{k+1} = n_{k+2} + n_{k+3} - p, \ n'_{k+2} = p, \ n'_{k+3} = n_{k+1} + n_{k+2} - p, \\
\text{where } \ p = \min( n_{k+1}, n_{k+3}).
\end{gathered}
\end{equation}

\item If $ a_{ij} = -1,  a_{ji} = -2 $, then $ d = 4 $ and
\begin{equation} \label{eq:valtrans3}
\begin{gathered}
n'_{k+1} = n_{k+2} + n_{k+3} + n_{k+4} - p_1, \ n'_{k+2} = 2p_1 - p_2, \\ \ n'_{k+3} = p_2 -p_1, \ n'_{k+4} = n_{k+1} + 2n_{k+2} + n_{k+3} - p_2 \\
\text{where } \ p_1 = \min( n_{k+1} + n_{k+2}, n_{k+1} + n_{k+4}, n_{k+3} + n_{k+4}), \\ 
p_2 = \min(n_{k+1} + 2n_{k+2},  n_{k+1} + 2n_{k+4},   n_{k+3} +  2n_{k+4}).
\end{gathered}
\end{equation}

\item If $ a_{ij} = -2, a_{ji} = -1 $ , then $ d = 4 $ and
\begin{equation} \label{eq:valtrans4}
\begin{gathered}
n'_{k+1} = n_{k+2} + 2 n_{k+3} + n_{k+4} - p_2, \ n'_{k+2} = p_2 - p_1, \\ \ n'_{k+3} = 2p_1 -p_2, \ n'_{k+4} = n_{k+1} + n_{k+2} + n_{k+3} - p_1 \\
\text{where } \ p_1 = \min( n_{k+1} + n_{k+2}, n_{k+1} + n_{k+4}, n_{k+3} + n_{k+4}), \\ 
p_2 = \min(2n_{k+1} + n_{k+2},  2n_{k+1} + n_{k+4},   2n_{k+3} +  n_{k+4}).
\end{gathered}
\end{equation}

\end{enumerate}

\end{Proposition}

\begin{proof}
If $ a_{ij} = 0 $ then the result holds with $ U = B(n_\bigdot) $.

Suppose that $ a_{ij} = a_{ji} = -1 $.  Let 
\begin{equation*} 
 U := \{ b_\bigdot \in B(n_\bigdot) : b^0_{k+1} + b^0_{k+3} \ne 0 \},
\end{equation*}
where $ b^0_l $ is the coefficient $t^{n_l} $ in $ b_l $.  
 
If $ b_\bigdot \in U $, then let $ b_\bigdot', \pi $ be determined from $ b_\bigdot $ by (\ref{eq:paramtrans2}). Since $ \pi = b_{k+1} + b_{k+3} $, we know that $ \val(\pi) = p $, as the leading terms of $b_{k+1}$ and $ b_{k+3} $ don't cancel.  In particular, the denominator $ \pi $  doesn't vanish.  Hence if $ b'_\bigdot $ is given by (\ref{eq:paramtrans2}), then by Proposition \ref{th:paramtrans}, $ \mathbf{y}_{\wi'}(b'_\bigdot) = \mathbf{y}_\wi(b_\bigdot) $.  Moreover, the valuation of the $ b'_l $ are given by (\ref{eq:valtrans2}), since $ \val(\pi) = p $.

The other cases follow similarly.
\end{proof}

Note that if $ \wi, \wi'$ are related by a braid move starting at position $ k $ and involving $i,j$, then $ \wi', \wi $ are related by a braid move starting at position $ k $ and involving $j,i$.  Moreover, let $ n'_\bigdot $  be the sequence of integers obtained from $ n_\bigdot $ by the formulas in Proposition \ref{th:valtrans}.  It is easy to see that $ n'_\bigdot $ is an $\wi'$-Lusztig datum of coweight $ \mu $.   It is also easy to see that $n_\bigdot $ is the sequence of integers obtained from $ n'_\bigdot $ by these formulas where we regard $ \wi', \wi $ as related by a braid move.

Now we transport our results from $ G(\K) $ to $\Gr $. 

\begin{Theorem} \label{th:localdense}
The intersection $ A^\wi(n_\bigdot) \cap A^{\wi'}(n'_\bigdot) $ is dense in $ A^\wi(n_\bigdot) $.
\end{Theorem}

\begin{proof}
Let $ U $ be the non-empty open subset of $ B(n_\bigdot) $ from Proposition \ref{th:valtrans}.  Since the map from $ B(n_\bigdot) $ to $ A^\wi(n_\bigdot) $ is surjective (Theorem \ref{th:BtoA}), the set $ Y = \{ [\mathbf{y}_\wi(b_\bigdot)] : b_\bigdot \in U \} $ is dense in $ A^\wi(n_\bigdot) $.  By Proposition \ref{th:paramtrans}, if $ L \in Y $, then $ L $ has a representative $ \mathbf{y}_{\wi'}(b'_\bigdot) $ for $ b'_\bigdot \in B(n'_\bigdot) $.  Hence by Theorem \ref{th:BtoA}, $ Y \subset A^{\wi'}(n'_\bigdot) $.  Hence the intersection is dense. 
\end{proof}

Note that the tropical Pl\"ucker relation (\ref{eq:A2trop}), (\ref{eq:B2trop1}), (\ref{eq:B2trop2}) at $ (w = w^\wi_k,i,j) $ only involves $ M_\gamma $ for $ \gamma $ an $\wi$ or $\wi'$-chamber weight.  This observation leads to the following result.

\begin{Proposition} \label{th:localplu}
Let $ L \in A^\wi(n_\bigdot) \cap A^{\wi'}(n'_\bigdot) $.  Then the collection $ \big( M_\gamma := \D_\gamma(L) \big)_{ \gamma \in \Gamma^\wi \cup \Gamma^{\wi'}} $ satisfies the tropical Pl\"ucker relation at $(w,i,j)$.
\end{Proposition}

\begin{proof}
If $ L \in A^\wi(n_\bigdot) \cap A^{\wi'}(n'_\bigdot) $, then we know $ \D_\gamma(L) $ for $ \gamma $ an $ \wi $ or $ \wi'$-chamber weight.  Since these are the only chamber weights which show up in the tropical Pl\"ucker relation, we just need to make a simple computation to check that the relation between $ n_\bigdot $ and $n'_\bigdot $ in Proposition \ref{th:valtrans} matches the tropical Pl\"ucker relation at $(w,i,j) $.

The case $ d= 2 $ is trivial because there is no tropical Pl\"ucker relation (in fact, in this case $ \Gamma^\wi = \Gamma^{\wi'} $).

Consider the case $ a_{ij} = a_{ji} = -1$.  Then by the length formula (\ref{eq:Mfromn}),
\begin{gather*}
n'_{k+2} = - M_{w \cdot \fund_i} - M_{w s_j s_i \cdot \fund_i} + M_{w s_j \cdot \fund_j} - \sum_{l \ne i, j} a_{li} M_{w \cdot \fund_l}, \\
n_{k+1} = -M_{w \cdot \fund_i} - M_{w s_i \cdot \fund_i} + M_{w \cdot \fund_j} - \sum_{l \ne i, j} a_{li} M_{w \cdot \fund_l}, \\
n_{k+3} = -M_{w s_i \cdot \fund_i} - M_{w s_j s_i \cdot \fund_i} + M_{w s_i s_j \cdot \fund_j} - \sum_{l \ne i, j} a_{li} M_{w \cdot \fund_l}.
\end{gather*}
By (\ref{eq:valtrans2}), $ n'_{k+2} = \min(n_{k+1}, n_{k+3}) $.  Substituting the above expressions into this equation gives 
\begin{align*}
-M_{w \cdot \fund_i} - M_{w s_j s_i \cdot \fund_i} + M_{w s_j \cdot \fund_j} = \min \big( &-M_{w \cdot \fund_i} - M_{w s_i \cdot \fund_i} + M_{w \cdot \fund_j}, \\
& -M_{w s_i \cdot \fund_i} - M_{w s_j s_i \cdot \fund_i} + M_{w s_i s_j \cdot \fund_j} \big)
\end{align*}
which is equivalent to the tropical Pl\"ucker relation (\ref{eq:A2trop}).

The other cases are similar.
\end{proof}

It is easy to see that the converse of this proposition holds, but we will not need this.
 
\subsection{Global picture}
Fix a coweight $ \mu \ge 0 $. Let $ \wi, \wi'$ be two reduced words related by a braid move involving $ i,j $, starting at position $ k $.  Let $ L \in X(\mu) $ and let $ n_\bigdot $, $ n_\bigdot' $ be the $ \wi,\wi'$-Lusztig data of $ P(L) $.  So $ L \in A^\wi(n_\bigdot) \cap A^{\wi'}(n'_\bigdot) $.  We say that $ L $ is $\wi, \wi'$-\textbf{generic} if $ n_\bigdot$ and $n'_\bigdot $ are related as in Proposition \ref{th:valtrans}.  By Proposition \ref{th:localplu}, if $ L $ is $ \wi, \wi' $-generic, then $ \D_\bigdot(L) $ satisfies  the tropical Pl\"ucker relation at $ (w^\wi_k,i,j) $.

We say that $ L \in X(\mu) $ is \textbf{generic} if $ L $ is $\wi, \wi'$-generic for every pair of reduced words $ \wi, \wi'$ related by a braid move.

If $ w \in W, i,j\in I $ are such that $ ws_i > w $ and $ ws_j > w $, then there exist a pair of reduced words $ \wi, \wi' $ related by a $ d$-move starting at position $ k = \operatorname{length}(w)$, involving $ i,j$  such that $w^\wi_k = w $. Visually, such $ w,i,j$ gives a 2-face of the permutahedron with lowest vertex $ w $ and this 2-face gives the transition between the reduced words $ \wi, \wi' $ (as in Figure \ref{fig:twopaths}).  Hence by Proposition \ref{th:localplu}, if $ L $ is generic, then $ \D_\bigdot(L) $ satisfies the tropical Pl\"ucker relation at $ (w, i, j) $.  Since $ (w,i,j) $ were arbitrary, $ \D_\bigdot(L) $ satisfies all the tropical Pl\"ucker relations.

\begin{Lemma} \label{th:gendense}
For any reduced word $ \wi $ and any $ \wi$-Lusztig datum $ n_\bigdot $, $ \{ L \in A^\wi(n_\bigdot) : L $ is generic$ \} $ is dense in $ A^\wi(n_\bigdot)$.  
\end{Lemma}

\begin{proof}
For any reduced word $\wj$ and any $\wj$-Lusztig datum $m_\bigdot$, define $ A^\wj_k(m_\bigdot) $ recursively by $A_0^\wj(m_\bigdot) := A^\wj(m_\bigdot) $ and for $ k > 0 $ by
\begin{equation*}
 A_k^\wj(m_\bigdot) := A_{k-1}^\wj(m_\bigdot) \cap \bigcap_{\wj'} A_{k-1}^{\wj'}(m'_\bigdot),
\end{equation*}
where  the intersection is over all reduced words $ \wj'$ which are related to $ \wj $ by a braid move and where $m'_\bigdot $ is the $\wj'$-Lusztig datum corresponding to $ m_\bigdot $ under Proposition \ref{th:valtrans}.  

We claim that for each $ k $, $A_{k+1}^\wj(m_\bigdot) $ is dense in $ A_k^\wj(m_\bigdot) $ and in $ A_k^{\wj'}(m'_\bigdot) $ whenever $ \wj $ and $ \wj'$ are related by a braid move and $m_\bigdot$ and $ m'_\bigdot$ are related as in Proposition \ref{th:valtrans}.  We proceed by induction. 

By Theorem \ref{th:localdense}, $ A^\wj(m_\bigdot) \cap A^{\wj'}(m'_\bigdot) $ is  dense in $ A^\wj(m_\bigdot) $ and $ A^{\wj'}(m'_\bigdot)$.  So $ A_1^\wj(m_\bigdot) $ is the intersection of subsets of $ A^{\wj}(m_\bigdot) $ which are  dense in $ A^{\wj}(m_\bigdot) $.  Moreover, these subsets are all constructible, hence $ A_1^\wj(m_\bigdot) $ is dense in $ A^{\wj}(m_\bigdot) $.  This also shows that  $A_1^{\wj}(m_\bigdot) $ is  dense in $ A^\wj(m_\bigdot) \cap A^{\wj'}(m'_\bigdot) $ and hence in $ A^{\wj'}(m'_\bigdot) $.  This establishes the base case. 

For the inductive step, let $ k > 0 $.  By the inductive hypothesis, $ A_k^\wj(m_\bigdot) $ and $ A_k^{\wj'}(m'_\bigdot) $ are each dense in each of $ A_{k-1}^\wj(m_\bigdot) $ and $ A_{k-1}^{\wj'}(m'_\bigdot) $.  Hence $A_k^\wj(m_\bigdot) \cap A_k^{\wj'}(m'_\bigdot) $ is dense in $A_{k-1}^\wj(m_\bigdot) \cap A_{k-1}^{\wj'}(m'_\bigdot) $ and so $ A_k^\wj(m_\bigdot) \cap A_k^{\wj'}(m'_\bigdot) $ is dense in $A^\wj_k(m_\bigdot) $ and in $ A^{\wj'}_k(m'_\bigdot)$ (since each of these are contained in $A_{k-1}^\wj(m_\bigdot) \cap A_{k-1}^{\wj'}(m'_\bigdot) $).  From here, the inductive step follows the same reasoning as the base case.

In these arguments, we are repeatedly using the fact that if $ U \subset V \subset X $ and if $ U $ is dense in $ X $, then $ U $ is dense in $ V $. 

Now, specialize to $ \wj = \wi, m_\bigdot = n_\bigdot $. Let $ \wj, \wj' $ be two reduced words which are related by a braid move and such that $ \wj $ is connected to $ \wi $ by less than $ k $ braid moves.  Suppose that $ L \in A_k^\wi(n_\bigdot) $, then by induction on $k$, we see that $ L $ is $\wj, \wj'$-generic.  Hence if $ k $ is larger than the largest number of braid moves needed to connect any two reduced words, then $ A_k^\wi(n_\bigdot) \subset \{ L \in A^\wi(n_\bigdot) : L \text{ is generic} \} $.  By a chain of dense inclusions, we see that $ A_k^\wi(n_\bigdot) $ is dense in $ A^\wi(n_\bigdot) $ and hence $\{L \in A^\wi(n_\bigdot) : L \text{ is generic}  \} $ is dense in $ A^\wi(n_\bigdot) $.
\end{proof}

\begin{proof}[Proof of Theorem \ref{th:BZcycle}]
Let $ \mu \ge 0 $ be a coweight and let $ M_\bigdot $ be a BZ datum of coweight $ (0,\mu) $.  Because of the action of $ \cwl$ is suffices to consider only this case.

Let $ \wi $ be a reduced word for $ w_0 $.  Let $n_\bigdot$ be the $\wi$-Lusztig datum corresponding to $ M_\bigdot $ under (\ref{eq:Mfromn}).  

If $ L \in A^\wi(n_\bigdot) $ is generic, then $ \D_\bigdot(L)  $ and $ M_\bigdot  $ both obey the tropical Pl\"ucker relations.  Moreover, they have the same values whenever $ \gamma $ is an $\wi$-chamber weight.  Suppose that $ \wi'$ is another reduced word, related to $ \wi $ by a $ d$-move involving $ i, j$ starting at position $ k $.  Then since both obey the tropical Pl\"ucker relation for $(w^\wi_k, i, j)$, we see that $ \D_\gamma(L) = M_\gamma $ whenever $ \gamma $ is a $ \wi'$-chamber weight.  Continuing this argument (and using the fact that any reduced word is connected to $ \wi $ by a sequence of braid moves), we see that $ \D_\gamma(L) = M_\gamma $ for all chamber weights $ \gamma$. So $ L \in A(M_\bigdot) $.

By Lemma \ref{th:gendense}, $ \{ L \in A^\wi(n_\bigdot) : L \text{ is generic} \} $ is dense in $ A^\wi(n_\bigdot) $ and by the above analysis, this set is contained in $ A(M_\bigdot) $, so we see that
\begin{equation*}
\overline{ \{ L \in A^\wi(n_\bigdot) : L \text{ is generic} \} } = \overline{A(M_\bigdot)} = \overline{ A^\wi(n_\bigdot)}.
\end{equation*}

By Theorem \ref{th:Ldirr}, $\overline{A^\wi(n_\bigdot)} $ is a component of $ \overline{X(\mu)} $, so $ \overline{A(M_\bigdot)} $ is a component.  So $ \overline{A(M_\bigdot)} $ is an MV cycle of coweight $ \mu $.

Conversely, if $ Z $ is a component of $ \overline{X(\mu)} $, then $ Z = \overline{A^\wi(n_\bigdot)} $ for some $ n_\bigdot $ by Theorem \ref{th:Ldirr}.  Let $ L \in A^\wi(n_\bigdot) $ be generic.  By the above analysis $ Z = \overline{A(M_\bigdot)} $.  Since $ L $ is generic,  $ ( M_\gamma = \D_\gamma(L)) $ satisfies the tropical Pl\"ucker relations.  Also $P(L) = P(M_\bigdot) $ is a pseudo-Weyl polytope, so $M_\bigdot $ satisfies the edge inequalities.  Finally, $ M_{\fund_i} = 0 $ for all $ i $, since $ L \in X(\mu) \subset S_e^0 $.  Hence $ M_\bigdot $ is a BZ datum of coweight $ (0, \mu) $.  So all MV cycles are of the desired form.
\end{proof}

\section{Minkowski sums of MV polytopes} \label{se:mink}
The MV polytopes for $ SL_3, Sp_4, SL_4 $ appeared without proof in \cite{jared2}.  Anderson expressed these MV polytopes by producing a finite list of prime MV polytopes such that every MV polytope was a Minkowski sum of these prime MV polytopes.  Moreover, he grouped these prime MV polytopes into ``clusters'', such that all Minkowski sum monomials of primes within a cluster were MV polytopes.

We will now show that for each $ G $, there exists such a finite set of prime MV polytopes.
Moreover, we will show how to find the primes and their groupings into clusters.

The proof of the following lemma is given at the end of appendix \ref{se:appendix}.

\begin{Lemma} \label{th:minksum}
If $P(M_\bigdot)$ and $ P(N_\bigdot) $ are two pseudo-Weyl polytopes, then so is their Min\-kowski sum $ P(M_\bigdot) + P(N_\bigdot) := \{ \alpha + \beta : \alpha \in P(M_\bigdot), \beta \in P(N_\bigdot) \}$.  Moreover $ P(M_\bigdot) + P(N_\bigdot) = P((M+N)_\bigdot) $.
\end{Lemma}

Combining Lemma \ref{th:minksum} with Theorem \ref{th:BZcycle}, we see that to understand Minkowski sums of MV polytopes, it is enough to understand sums of BZ data.  In what follows, we will identify MV polytopes with their BZ data, so we will use $ \mathcal{P} $ to denote the set of BZ data.

If $ M_\bigdot$, $ N_\bigdot $ are BZ data, then $ (M+N)_\bigdot $ is not necessarily a BZ datum.  We will now see how to divide the set of BZ data into regions, within which we can add BZ data.  In this section, BZ datum always means a BZ datum of coweight $ (0, \cdot) $, so $ M_{\fund_i} = 0 $ for all $ i $.

\subsection{Prime BZ data}
If $ A = \min(B,C, D) $ is a $(\min, +)$ equation, then a \textbf{min-choice} for this equation is a choice of $ B, C, $ or $ D $.  Corresponding to such a choice, we get a system of linear equations and inequalities of the form
\begin{equation*}
A \ge B, \quad A = C, \quad A \ge D.
\end{equation*}
Note that if $ A,B,C,D $ is a solution to the original $(\min, +)$ equation, then it satisfies the system corresponding to at least one of the three possible min-choices.  In fact, the (non-disjoint) union of all solutions to the three systems is the set of solutions to the original equation.

A \textbf{BZ-choice} is a collection of min-choices, one for each tropical Pl\"ucker relation.  Note that there are $ 2^{\#H} 9^{\#O} $ possible BZ choices, where $ \#H$ and $ \#O $ are the number of hexagons and octagons in the permutahedron for $ G $.

  Let $ \Sigma $ denote the set of BZ-choices.  If $ \sigma \in \Sigma $, then let $ \mathcal{P}_\sigma $ denote the set of BZ data which satisfy the systems corresponding to each min-choice in $ \sigma $.  Note that $\mathcal{P} = \cup_{\sigma \in \Sigma} \mathcal{P}_\sigma $, but that this union is not disjoint.  

Each $\mathcal{P}_\sigma $ is the set of lattice points of a rational polyhedral cone in $ \R^\Gamma $ --- namely the cone defined by all the linear equations and inequalities coming from min-choices in $ \sigma $, by the edge inequalities, and by the equations $ M_{\fund_i} = 0$ for all $ i $.  Moreover, this cone lies in $ \R_{\le 0}^\Gamma $, since the edge inequalities imply that $ M_\gamma \le 0 $ for all $ \gamma$.  Hence it is a proper cone.  Since $ \mathcal{P}_\sigma $ is the set of lattice points of a cone, if $ M_\bigdot, N_\bigdot \in \mathcal{P}_\sigma $, then $ (M+N)_\bigdot \in \mathcal{P}_\sigma $.  So $ \mathcal{P}_\sigma $ forms a monoid.    By Gordan's Lemma (see \cite[section 3]{E}), the monoid is finitely generated by a unique minimal set of generators which we call the $\sigma$-\textbf{prime BZ data}.

\subsection{Prime MV polytopes}
A $\sigma$-\textbf{prime MV polytope} is an MV polytope corresponding to a $\sigma$-prime BZ datum.  Each set of $\sigma$-prime MV polytopes is called a \textbf{cluster} of prime MV polytopes.  Thus the clusters are indexed by the set $ \Sigma $ of possible BZ choices.  There are finitely many clusters and finitely many prime MV polytopes in each cluster.

Combining the above observations gives the following result.

\begin{Theorem} \label{th:MVprime}
Every MV polytope is the Minkowski sum of prime MV polytopes.  More specifically, every MV polytope is sum of $\sigma$-prime MV polytopes for some $\sigma$.  Moreover, the sum of $\sigma$-prime MV polytopes is always an MV polytope.
\end{Theorem}

As noted above, this result was first observed in low rank cases by Anderson \cite{jared2}.  In \cite{jaredmisha}, Anderson-Kogan argued that the existence of this canonical set of generators for the set of MV polytopes is related to the cluster algebras of Berenstein-Fomin-Zelevinsky \cite{BFZ2}.  This connection is an interesting direction for future research.  See \cite{jaredmisha2} for recent results in this direction.

Note that not all the cones corresponding to different BZ choices are of the same dimension.  In low rank examples, we have observed that there are some (very few) cones of dimension $ m $ and all the rest of the cones are subcones of these maximal cones.  Moreover, not all the maximal cones are isomorphic.  In general, it seems to be an interesting problem to understand the structure of these cones.

The analysis of the case of $ SL_3 $ is easy and was carried out in section \ref{se:MVpoly}.  The tropical Pl\"ucker relation (\ref{eq:plusl3}) gives two BZ choices, each of which gives a maximal cone.  The two maximal cones give the two kinds of $ SL_3 $ MV polytopes shown in figure \ref{fig:sl3}.

\begin{figure}
\begin{center}
\epsfig{file=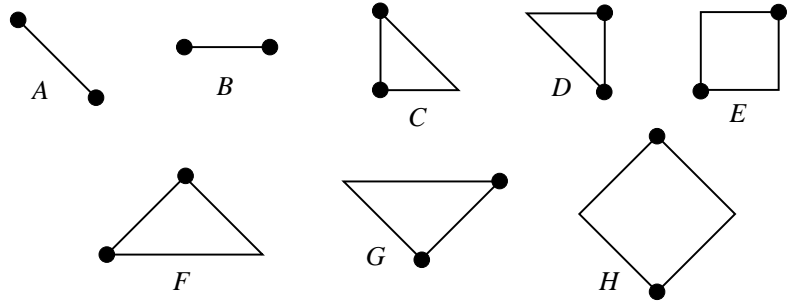,height=4cm}
\caption{The prime $Sp_4$ MV polytopes.}
\label{fig:primesp4}
\end{center}
\end{figure}

For $ Sp_4 $, only 4 of the 9 cones are maximal and these lead to the
four possible types of $Sp_4 $ MV polytopes shown in figure \ref{fig:sp4}.   Two of the maximal cones have 4 generators (these are simplicial) and the other two cones have 5 generators.  Figure \ref{fig:primesp4} shows the 8 prime $Sp_4 $ MV polytopes labelled by the letters $ A, \dots, H $.  The maximal clusters of prime $Sp_4 $ MV polytopes (written in the order that they correspond to the types in figure \ref{fig:sp4}) are $ BEFGH $, $CEFH $, $ACDEH $, and $ DEGH $.  These polytopes and clusters first appeared in \cite{jared2}.  Note that $ E $ and $ H $ appear in each maximal cluster.  This is because they are the Weyl polytopes for the fundamental weights.

For $ SL_4 $,  there are $2^8 = 256$ cones, 13 of which are maximal.  Of these 13 maximal cones, 12 have 6 generators and are simplicial, while 1 has 7 generators.  There are a total of 12 prime MV polytopes in this case.

\section{Relation to the canonical basis}
\subsection{Lusztig data} \label{se:lusdat}
Let $ \wi $ be a reduced word for $ w_0 $.  If $ P $ is an MV polytope, then we can extract its $\wi$-Lusztig datum.  This is invariant under the $ \cwl $ action and so the $ \wi$-Lusztig datum of a stable MV polytope is well-defined.  

\begin{Theorem} \label{th:lusdatbij}
Taking $ \wi$-Lusztig datum gives a bijection $ \psi_\wi : \mathcal{P} \rightarrow \N^m $.
\end{Theorem}
\begin{proof}
By Theorem \ref{th:Ldirr}, an inverse is given by $ n_\bigdot \mapsto \Phi \big(\overline{A^\wi(n_\bigdot)} \big)$.
\end{proof}

Suppose that $ \wi $ and $ \wi'$ are two reduced words for $ w_0 $.  Then the transition map $ \psi_{\wi'} \circ \psi_{\wi}^{-1} : \N^m \rightarrow \N^m $ is a bijection.  When $ \wi, \wi' $ are related by a braid move, then Theorem \ref{th:localdense} shows that this bijection is given by the $(\min, +) $ equations in Proposition \ref{th:valtrans}.  As shown in the proof of Proposition \ref{th:localplu}, these bijections are equivalent to the tropical Pl\"ucker relations.  In practice, when trying to construct MV polytopes, these bijections are often easier to work with than the tropical Pl\"ucker relations.  However, we have started with the tropical Pl\"ucker relations because they are more naturally motivated and are often better for theoretical purposes.

\subsection{The canonical basis}
Recall that $ G^\vee $ is the group with root datum dual to that of $ G $.  In particular, the weight lattice of $ G^\vee $ is $ \cwl $.  Let $ \mathcal{B} $ denote Lusztig's canonical basis for $ U^\vee_+$, the upper triangular part of the quantized universal enveloping algebra of $ G^\vee$.  Lusztig showed that a choice of reduced word $ \wi $ for $ w_0 $ gives rise to a bijection $ \phi_{\wi} : \mathcal{B} \rightarrow \N^m $ (see \cite[section 2]{Lbook} or \cite[Proposition 4.2]{BZtpm} for more details). Following Berenstein-Zelevinsky, we call $\phi_{\wi}(b) $ the $\wi$-\textbf{Lusztig datum} of $ b $.

Moreover, Lusztig \cite[section 2.1]{Lbook} and \cite[section 12.5]{Lproj} showed that the transition map $\phi_{\wi'} \circ \phi_{\wi}^{-1} $ matches the bijection in Proposition \ref{th:valtrans}, whenever $ \wi, \wi' $ are related by a braid move.  In fact, this was the original source of these bijections.  Since any two reduced words are connected by a sequence of braid moves, we see that $ \phi_{\wi'} \circ \phi_{\wi}^{-1} = \psi_{\wi'} \circ \psi_{\wi}^{-1} $ for all reduced words $ \wi, \wi'$ and we immediately obtain the following result.

\begin{Theorem} \label{th:can}
There is a coweight preserving bijection $ b \mapsto P(b) $ between the canonical basis $ \mathcal{B} $ and the set $\mathcal{P} $of stable MV polytopes.  Under this bijection, the $\wi$-Lusztig datum of $ b $ equals the $\wi$-Lusztig datum of $P(b)$.
\end{Theorem}
In other words to find the $\wi$-Lusztig datum of $ b $, we can just look at the lengths of the edges in $P(b) $ along the path determined by $ \wi $.

In fact, Lusztig noticed in \cite{Lbookbook} that the transition map $ \phi_{\wi'} \circ \phi_{\wi}^{-1} : \Z^m \rightarrow \Z^m $ was the tropicalization of the transition map $ y_{\wi'}^{-1} \circ y_{\wi} : \R^m \rightarrow \R^m $ between the parametrizations of $ N $ (these are given in our Proposition \ref{th:paramtrans}).  Lusztig and Berenstein-Zelevinsky further explored this relationship in \cite{Lfunnypaper} and \cite{BZtpm} respectively.  In the latter paper, which served as a primary motivation for this work, Berenstein-Zelevinsky invented the notion of BZ data.  More specifically, combining Theorem 4.3 in \cite{BZschub} and Example 5.4 in \cite{BZtpm}, they showed that there is a bijection between the set of BZ data of coweight $ (0,\cdot) $ and the canonical basis.  This motivated us to look for a bijection between BZ data and MV cycles.  

\section{Finite dimensional representations}
As stated in the introduction, one of the main purposes of studying MV cycles and polytopes is to use them to understand the combinatorics of finite-dimensional representations of $ G^\vee $.

\subsection{Indexing the MV basis}
If $ \lambda \in \cwl^+$, then let $\Gr^\lambda := \overline{t^\lambda G(\mathcal{O})} $.  It is known that $ \Gr^\lambda $ is a finite-dimensional projective variety.  The geometric Satake isomorphism gives an isomorphism of the intersection homology $ \mathrm{IH}(\Gr^\lambda) $ with the finite-dimensional representation $ V_\lambda $.  So intersection homology cycles in $ \Gr^\lambda $ give us elements of $ V_\lambda $.  Mirkovic-Vilonen \cite{MV} proved that the MV cycles were good cycles to consider.

\begin{Theorem} \label{th:MVbasis}
Under the isomorphism $ \mathrm{IH}(\Gr^\lambda) \overset{\sim}{\rightarrow} V_\lambda $, the MV cycles of coweight $ (\mu, \lambda) $ which lie in $ \Gr^\lambda $ give a basis (the \textbf{MV basis}) for the weight space  $ V_\lambda(\mu) $.
\end{Theorem}

The moment map image of $ \Gr^\lambda $ is the Weyl polytope $ W_\lambda $.  Hence if $ A $ is an MV cycle lying in $\Gr^\lambda $, then its moment polytope $ \Phi(A)$ lies in $ W_\lambda $.  Anderson showed that the converse holds.

\begin{Lemma}[{\cite[Proposition 7]{jared2}}] \label{th:jaredcontain}
Let $ A $ be an MV cycle of coweight $ (\mu, \lambda) $.  Then $ A \subset \Gr^\lambda $ if and only if $ \Phi(A) \subset W_\lambda $.
\end{Lemma}

Recall that $ W_\lambda $ is a pseudo-Weyl polytope corresponding to the collection $ N_{w \cdot \fund_i} = \langle w_0 \cdot \lambda, \fund_i \rangle $.  Hence combining Lemma \ref{th:jaredcontain} and our Theorem \ref{th:BZcycle} immediately gives the following result.

\begin{Theorem} \label{th:MVwtm}
The MV basis for $ V_\lambda(\mu) $ is indexed by the set of BZ data $ M_\bigdot $ of coweight $ (\mu, \lambda) $ such that 
\begin{equation} \label{eq:cond}
M_{w \cdot \fund_i} \ge \langle w_0 \cdot \lambda, \fund_i \rangle,
\end{equation}
for all $ i \in I $ and $ w \in W $.
In particular, counting such BZ data gives a formula for the weight multiplicity.
\end{Theorem}

\subsection{Indexing the canonical basis}
On the other hand, we can also consider the canonical basis (specialized at $ q=1$) for a representation $ V_\lambda $.  One way to describe this basis is to consider the map $ \eta_\lambda : U(\mathfrak{n}^\vee) \rightarrow V_\lambda $ which is given by acting on the low weight vector.  Let $ \mathcal{B}(\lambda) $ denote the subset of $\mathcal{B} \subset U(\mathfrak{n}^\vee) $ which is not sent to 0 by this map.  Lusztig \cite[section 8]{Lbook} proved that $ \eta_\lambda $ maps $ \mathcal{B}(\lambda) $ bijectively onto a basis for $ V_\lambda $, which is called the \textbf{canonical basis} for $ V_\lambda $.  Moreover, this basis is compatible with weight spaces.

Lusztig also characterized $ \mathcal{B}(\lambda) $ in terms of Lusztig data.  In particular, only the last component $\phi_\wi(b)_m $ is relevant.
\begin{Theorem}[{\cite[section 8]{Lbook}, \cite[Corollary 3.4]{BZtpm}}] \label{th:BZcriteria} Let $ b \in \mathcal{B} $.  Then
$ b \in B(\lambda) $ if and only if 
\begin{equation*}
 \phi_\wi(b)_m \le -\langle w_0 \cdot \lambda, \alpha_{i_m} \rangle
\end{equation*}
for all reduced words $ \wi$.
\end{Theorem}

Using the bijection between MV polytopes and the canonical basis (Theorem \ref{th:can}) and the description of MV polytopes by BZ data (Theorem \ref{th:BZcycle}), we immediately obtain the following result which is a version of Theorem 5.16 from \cite{BZtpm}.

\begin{Theorem} \label{th:canwtm}
The canonical basis for $ V_\lambda(\mu) $ is indexed by the set of BZ data $ M_\bigdot $ of coweight $ (\mu, \lambda) $ such that 
\begin{equation*}
M_{w_0 s_i \cdot \fund_i} \ge \langle w_0 \cdot \lambda,\fund_i \rangle.
\end{equation*}
In particular, counting such BZ data gives a formula for the weight multiplicity.
\end{Theorem}

It interesting to compare Theorems \ref{th:MVwtm} and \ref{th:canwtm}.  The condition on the BZ data in Theorem \ref{th:canwtm} is appears weaker since we only impose (\ref{eq:cond}) for $ w = w_0 s_i $ for some $ i $.  In other words, in Theorem \ref{th:MVwtm} we demand that all vertices $ \mu_\bigdot $ of the polytope lie in $W_\lambda $, whereas in Theorem \ref{th:canwtm} we only require that $ \mu_{w_0 s_i} \in W_\lambda $ for all $ i $.  Hence the set of BZ data in Theorem \ref{th:canwtm} is \textit{a priori} bigger than that in Theorem \ref{th:MVwtm}.  However, the two sets of BZ data index bases for the same finite dimensional vector space, hence they must be the same set.  

In particular all the inequalities (\ref{eq:cond}) for $ w \ne w_0 s_i $ are redundant.  It would be interesting to find a direct combinatorial proof of this fact.  Such a proof seems to require a good understanding of the combinatorics of the tropical Pl\"ucker relations.  We have been able to find such a proof for $ SL_n $ and $ SO_{2n} $ but not in general. 

\subsection{Tensor product multiplicities}
Anderson also extended Theorem \ref{th:MVbasis} to show that MV cycles give a basis for tensor product multiplicity spaces.  Using Lemma \ref{th:jaredcontain}, he proved the following tensor product multiplicity formula.

\begin{Theorem}[{\cite[Theorem 1]{jared2}}] \label{th:jaredtpm}
Let $ \lambda, \mu, \nu \in \cwl^+ $.
The tensor product multiplicity $c_{\lambda\mu}^\nu $ of $ V_\nu $ inside $ V_\lambda \otimes V_\mu $ is equal to the number of MV polytopes $ P $ such that
\begin{enumerate}
\item $ P $ has coweight $ (\nu- \mu, \lambda) $,
\item $P $ is contained in $W_\lambda $,
\item $P $ is contained in $ - W_\mu + \nu $.
\end{enumerate}
\end{Theorem}

Combining Theorem \ref{th:jaredtpm} with Theorem \ref{th:BZcycle}, we immediately obtain the following result.
\begin{Theorem} \label{th:BZtpm}
The multiplicity $ c_{\lambda\mu}^\nu $ equals the number of BZ data of coweight $ (\nu - \mu, \lambda) $ such that 
\begin{enumerate}
\item $ M_\gamma \ge \langle w_0 \cdot \lambda, \fund_i \rangle$ for all $ i $ and for all chamber weight $ \gamma $ of level $ i $, 
\item $ M_\gamma \ge \langle \nu, \gamma \rangle  - \langle \mu, \fund_i \rangle $ for all $ i $ and for all chamber weights $ \gamma $ of level $ i $.
\end{enumerate}
\end{Theorem}

As with weight multiplicity, there is also a tensor product multiplicity formula coming from the canonical basis.  This is given as Theorem 5.16 in \cite{BZtpm}.  It can be obtained from the above theorem by considering only those chamber weights of the form $ w_0 s_i \cdot \fund_i $ and $ s_i \cdot \fund_i $ in (i) and (ii) respectively.  The relationship between these two tensor product multiplicity formulas is the same as the previously discussed relationship between the corresponding weight multiplicity formulas (Theorems \ref{th:MVwtm} and \ref{th:canwtm}).

\section{$SL_n$ comparison} \label{se:typeA}

We now examine our constructions in greater detail when $ G = SL_n$.  The main goal of this section is to connect our work with that of Anderson-Kogan \cite{jaredmisha}.  We also hope that the ideas presented here will help the reader get a better feel for our main results.

\subsection{Lattices}
Let $ U $ be a vector space over $\C $.  A \textbf{lattice} in $ U \otimes \K $ is a free $ \oo $-submodule $L \subset U \otimes \K $ such that $ \Span_\K(L) = U \otimes \K $. Let $ L_0 := U \otimes \mathcal{O} $ denote the \textbf{standard lattice} in $ U \otimes \K $.  The \textbf{relative dimension} of a lattice $ L $ in $ \mathcal{U} \otimes \K $ is defined to be $ \dim_\C(L/L \cap L_0 ) - \dim_\C(L_0 / L \cap L_0 )$ and is denoted $ \rdim (L) $.

If $ G $ is a reductive group and $ V_\lambda $ is a representation of $ G $, then there is a map
\begin{equation} \label{eq:grtolat}
\begin{aligned}
\psi_\lambda : \Gr &\rightarrow \{ \text{lattices in } V_\lambda \otimes \K \} \\
[g] &\mapsto g^{-1} \cdot V_\lambda \otimes \oo.
\end{aligned}
\end{equation}

For any $ G $, this gives an embedding of $ G $ into $ \prod_{\lambda} \{ \text{lattices in } V_\lambda \otimes \K \} $.  The image of this embedding will be those systems of lattices which are compatible with morphisms $ V_\lambda \otimes V_\mu \rightarrow V_\nu $ (see section 10.3 in \cite{FM} for more details).  We can use this embedding to express our functions $ \D_\gamma $.

\begin{Proposition} \label{th:Dgamlatgen}
Let $ \gamma $ be a chamber weight of level $ i $.  Then
\begin{equation*}
 \D_\gamma(L) = \rdim \big(\psi_{\fund_i}(L) \cap (V_{\fund_i}(\gamma) \otimes \K) \big).
 \end{equation*}
\end{Proposition}

\begin{proof}
Note that if $ R $ is a lattice in $ V_{\fund_i} $, then since $ V_{\fund_i}(\gamma) $ is one dimensional,
\begin{equation*}
\rdim \big(R \cap (V_{\fund_i}(\gamma) \otimes \K) \big) =  -\min \val(a) 
\end{equation*}
where the $\min$ is taken over all $ a $ such that $ a v_\gamma \in R $.  

Hence, in our case the $\min $ is taken over all $ a $ such that 
\begin{equation*}
 a v_\gamma \in g^{-1} \cdot V_{\fund_i} \otimes \mathcal{O} \ \Leftrightarrow \ a g \cdot v_\gamma \in V_{\fund_i} 
 \otimes \mathcal{O}.
 \end{equation*} 
 From here the result follows from the definition of $ \val $.
\end{proof}

\subsection{Lattices for $ SL_n$}
From now on, we specialize to the case $ G = SL_n $ where an easier picture is available.  The following result is due to Lusztig.

\begin{Theorem}[\cite{L}] \label{th:latmodel}
In the case of the standard representation $ V_{\fund_1} $ of $ SL_n$, the map $ \psi_{\fund_1} $ gives an isomorphism $ \Gr \rightarrow \Grl := \{ \text{lattices in } \K^n \text{ of relative dimension } 0 \} $.
\end{Theorem}

If $ U $ is a vector space over $ \C $ and $ L $ is a lattice in $ U \otimes \K $, then $ \Lambda^i L := \{ v_1 \wedge \cdots \wedge v_k : v_1, \dots, v_i \in L \} $ is a lattice in $ \Lambda^i U \otimes \K $.  Since the $ i$th  fundamental representation of $ SL_n $ is $ \Lambda^i \C^n $, we see that if $ L \in \Gr $, then
\begin{equation} \label{eq:wedge}
 \psi_{\fund_i}(L) = \Lambda^i \big( \psi_{\fund_1}(L) \big).
\end{equation}
So from $ \psi_{\fund_1}(L) $ we can recover $ \psi_{\fund_i}(L) $ for all $ i $, and then from there $ \psi_\lambda(L) $ for all $ \lambda $.

Let $ \{e_1, \dots, e_n \} $ denote the usual basis for $ \C^n $.  Recall that for each $ \mu = (\mu_1, \dots, \mu_n) \in \cwl $, we defined an element $ t^\mu \in \Gr $.  Under the above isomorphism, this element goes over to the lattice $ L_\mu := \Span_{\mathcal{O}}(t^{-\mu_1} e_1, \dots, t^{-\mu_n} e_n ) $.

Recall that the set of chamber weights $ \Gamma $ can be identified with the set of proper subsets of $ \{ 1, \dots, n \} $.  For any $ \gamma \in \Gamma $, we can consider the subspace $ U_\gamma := \Span\{ e_i : i \in \gamma \} $ of $ \C^n $.  We get a corresponding subspace $ U_\gamma \otimes \K $ of $ \K^n $.  In what follows we will abuse notation and write $ U_\gamma $ for $ U_\gamma \otimes \K $.  

\begin{Proposition} \label{th:dgam}
Under the isomorphism $ \Gr \rightarrow \Grl $, the function $ \D_\gamma $ becomes the function
\begin{equation*}
L \mapsto \rdim(L \cap U_\gamma)
\end{equation*}
\end{Proposition}

The proof of this proposition follows from Proposition \ref{th:Dgamlatgen}, equation (\ref{eq:wedge}), and the following lemma.

\begin{Lemma} \label{th:wedge}
If $ U $ is a vector space over $\C $ of dimension $ k $ and $ L $ is a lattice in $ U \otimes \K $, then $ \rdim(L) = \rdim(\Lambda^k L) $.
\end{Lemma}

\begin{proof}
Fix a basis $ \{u_1, \dots, u_k \} $ for $ U $ over $ \C $.  

Suppose that there exist $ r_1, \dots, r_k \in \Z$ such that $ L = \Span_\mathcal{O}(t^{-r_1} u_1, \dots, t^{-r_k} u_k) $.  Then it is easy to see that $ \rdim (L) = r_1 + \cdots + r_k $ and
\begin{gather*} 
 \rdim(\Lambda^k L) = \rdim( \Span_\mathcal{O}(t^{-r_1 - \dots - r_k} u_1 \wedge \dots \wedge u_k)) = r_1 + \cdots + r_k .  
\end{gather*}

If $ L $ is a lattice, then there exists $ g \in GL_U(\mathcal{O}) $ such that $ g \cdot L $ is of the above form.  Hence it suffices to show that $ \rdim(L) $ and $ \rdim(\Lambda^k L) $ are invariant under $ g $. 

Note that $ g \cdot L_0 = L_0 $, so $ g \cdot (L \cap L_0) = (g \cdot L) \cap L_0 $.  Hence
\begin{equation*}
\dim(L/L \cap L_0) = \dim(g \cdot L / (g \cdot L) \cap L_0) \text{ and } \dim(L_0/L \cap L_0) = \dim(L_0/ (g\cdot L) \cap L_0).  
\end{equation*}
So $ \rdim(L) $ is invariant under $ g $.  Similarly, $\rdim(\Lambda^k L) $ is invariant under $ g $. 
\end{proof}

\subsection{Kostant pictures}
The positive roots of $ SL_n $ are indexed by the set of pairs $ \{ (a,b) : 1 \le a < b \le n \} $ with the positive root corresponding to $ (a,b) $ having a $ 1 $ in the $a$th slot, a $ -1$ in the $b$th slot and $0$s elsewhere.

Following Anderson-Kogan \cite{jaredmisha}, we define a \textbf{Kostant picture} to be an assignment of non-negative integer to each positive root of $ SL_n $.  They viewed a Kostant picture as a collection of loops around the Dynkin diagram for $ SL_n $, but we will think of it more formally as an element $ \mathbf{p} = \big( p_{(a,b)} \big)_{a < b} \in \N^{\Delta_+}$.  To reconstruct the ``picture'', one should draw $ p_{(a,b)} $ loops with left edge at column $ a $ and right edge at column $ b $.

In \cite{jaredmisha}, Anderson-Kogan produced a bijection from the set of Kostant pictures to the set of stable MV cycles and polytopes for $ SL_n $.  Using $\wi $-Lusztig data, we constructed such a bijection for any reduced word $ \wi $ for $w_0$ (Theorems \ref{th:Ldirr} and \ref{th:lusdatbij}).  We will show that their bijection is our bijection for the reduced word $ \wi := (1, \dots, n-1, 1, \dots, n-2, \dots, 1,2,1) $.  

Recall that any reduced word induces a total order on the positive roots.    The reduced word $ \wi $ induces the order 
\begin{equation*} 
\beta_1^\wi = (1,2), (1, 3), \dots, (1,n), (2,3), \dots, (2,n), \dots, (n-1, n) = \beta_m^\wi .
\end{equation*}

Recall that any reduced word  $\wi $ has its associated set of $\wi$-chamber weights, denoted $ \Gamma^\wi $.  For this choice of $ \wi $, we see that $ \Gamma^\wi := \{ [a \cdots b] : a < b \} $ where $ [a \cdots b] $ denotes the chamber weight $ \{a, a+1, \dots, b \} $. 

We now recall some more notation from \cite{jaredmisha}.  If $ L $ is a lattice, then $ \delta_i(L) $ denotes the maximum value of $ j $ such that $ t^{-j} e_i \in L $.  So $ L_{\delta(L)} \subset L $.  Next, Anderson-Kogan define $ \dim_0(L) := \dim_\C(L / L_{\delta(L)}) $.

The following lemma concerning these functions will be very convenient for us.

\begin{Lemma} \label{th:dim0}
\begin{equation*}
\dim_0(L \cap U_\gamma) + \sum_{i \in \gamma} \delta_i(L) = \rdim(L \cap U_\gamma)
\end{equation*}
\end{Lemma}

\begin{proof}
If $ A \subset B \subset C $ is a tower of vector spaces over $ \C $, then $ \dim(C/A) = \dim(C/B) + \dim(B/A) $.  In our case, we have the three towers $ L_0 \cap L_{\delta(L)} \subset L_0 \cap L \subset L_0 $, $ L_0 \cap L_{\delta(L)} \subset L_{\delta(L)} \subset L $, and $ L_0 \cap L_{\delta(L)} \subset L_0 \cap L \subset L $.  Applying this tower theorem to the intersection of these three towers with $ U_\gamma $ and adding the equations implies the desired result.
\end{proof}

If $ (a,b) $ is a positive root of $ SL_n$ (which they think of as a loop), Anderson-Kogan defined a function $ n_{(a,b)} : \Grl \rightarrow \Z $ by the formula
\begin{equation*}
n_{(a,b)}(L) := \dim_0(L \cap U_{[a \cdots b]}) - \dim_0(L \cap U_{[a \cdots b-1]}) - \dim_0(L \cap U_{[a+1 \cdots b]}) + \dim_0(L \cap U_{[a+1 \cdots b-1]}).
\end{equation*}

\begin{Proposition}
Let $ L \in \Gr $ and let $ n_\bigdot  $ denote the $ \wi $-Lusztig datum of $ P(L) $.  Let $ L $ also denote the image of $ L $ in $ \Grl $.  Let $ (a,b) $ be a positive root and let $ k $ be such that $ \beta_k^\wi = (a,b) $.  Then
\begin{equation*}
n_{(a,b)}(L) = n_k.
\end{equation*}
\end{Proposition}

\begin{proof}
Let $ M_\gamma = \D_\gamma(L) $.  Examining the reduced word $ \wi $ and the system (\ref{eq:Mfromn}) which converts between $\wi$-Lusztig datum and hyperplane datum for $\wi$-chamber weights, we see that 
\begin{equation*}
n_k = - M_{[a+1 \cdots b]} - M_{[a \cdots b-1]} + M_{[a+1 \cdots b-1]} + M_{[a \cdots b]}.
\end{equation*}
 
So it suffices to prove that
\begin{equation*}
n_{(a,b)}(L) = -\D_{[a+1 \cdots b]}(L) - \D_{[a \cdots b-1]}(L) + \D_{[a+1 \cdots b-1]}(L) + \D_{[a \cdots b]}(L)
\end{equation*}
for all $ L \in \Gr $.
But this follows immediately from Proposition \ref{th:dgam} and Lemma \ref{th:dim0}.
\end{proof}

If $ \mathbf{p} = (p_{(a,b)}) $ is a Kostant picture, then Anderson-Kogan defined the notion of a lattice \textbf{weakly compatible} to $ \mathbf{p}$.  They showed that a lattice $L $ was weakly compatible to $ \mathbf{p} $ iff $ n_{(a,b)}(L) = p_{(a,b)} $ for all $ (a,b) $ \cite[Prop 4.1]{jaredmisha}.  Following a slight modification of their notation, we let $ M(\mathbf{p}) $ denote the set of lattices $ L $ weakly compatible to $ \mathbf{p} $ such that $ \rdim(L \cap U_{[1 \cdots b]}) = 0 $ for all $b$.  (Actually, Anderson-Kogan considered ``$ M(\mathbf{p}, \lambda) $'' which consisted of lattices $ L $ weakly compatible to $ \mathbf{p} $ and satisfying the condition $ \rdim(L \cap U_{[a \cdots n]}) = \langle \lambda, [a \cdots n] \rangle $ for all $ a $). 

From the above proposition, we immediately see how these lattices fit into our results. 

\begin{Corollary}
Let $ n_\bigdot \in \N^m $.  Define $ \mathbf{p} $ by $ p_{(a,b)} = n_k $ where $ k $ is such that $ \beta^\wi_k  = (a,b)$.  The isomorphism $ \Gr \rightarrow \Grl $ takes $ A^\wi(n_\bigdot) $ onto $ M(\mathbf{p}) $.
\end{Corollary}

Anderson-Kogan proved that the closure of $ M(\mathbf{p}) $ was an MV cycle.  Our Theorem \ref{th:Ldirr} shows that $\overline{A^\wi(n_\bigdot)} $ is an MV cycle and thus provides an alternate proof of this result.   Thus we have established the following.

\begin{Theorem} \label{th:AKbij}
The Anderson-Kogan bijection $ \N^{\Delta_+} \rightarrow \mathcal{M} $ is the same as our bijection $ \N^m \rightarrow \mathcal{M} $ for the particular choice of reduced word above.
\end{Theorem}

\subsection{Strongly compatible lattices}

In order to understand the MV polytope $\Phi \big(\overline{M(\mathbf{p})} \big) $ associated to $ \overline{M(\mathbf{p})} $, Anderson-Kogan introduced the notion of a lattice in $ M(\mathbf{p}) $ being \textbf{strongly compatible} to $ \mathbf{p}$.  We follow their notation and write $ M^\ddag(\mathbf{p}) $ for the set of lattices in $ M(\mathbf{p}) $ which are strongly compatible to $\mathbf{p}$.  

They showed that $ M^\ddag(\mathbf{p}) $ was dense in $ M(\mathbf{p}) $ \cite[Prop 4.3,4.5]{jaredmisha} and that a lattice $ L $ was strongly compatible to $ \mathbf{p} $ if and only if $ P(L) = \Phi \big( \overline{M(\mathbf{p})} \big) $ \cite[Theorem E]{jaredmisha}.  Anderson-Kogan also observed that $ M^\ddag(\mathbf{p}) $ was a GGMS stratum and this served as one of our primary motivations for the use of GGMS strata in this work.

We showed that for any $ \wi$-Lusztig datum $ n_\bigdot $ of coweight $ \mu$, there exists a corresponding BZ datum $ M_\bigdot $ of coweight $ (0, \mu) $ (Theorem \ref{th:lusdatbij}).  Moreover, the corresponding GGMS stratum $ A(M_\bigdot) $ is dense in $ A^\wi(n_\bigdot) $ and the corresponding pseudo-Weyl polytope $ P(M_\bigdot) $ equals the MV polytope $ \Phi \big( \overline{A^\wi(n_\bigdot)} \big) $.  For these results see Lemma \ref{th:gendense} and the proof of Theorem \ref{th:BZcycle}.  

Comparing our results to the Anderson-Kogan results, we immediately have the following.

\begin{Theorem}
Let $ n_\bigdot \in \N^m $ and let $ \mathbf{p} $, $ M_\bigdot $ be related to $n_\bigdot $ as above.  Then the isomorphism $ \Gr \rightarrow \Grl $ takes $ A(M_\bigdot) $ onto $ M^\ddag(\mathbf{p}) $.
\end{Theorem}

\subsection{Collapse algorithm}

Fix $ n_\bigdot $ and the corresponding $ \mathbf{p} $ and $ M_\bigdot$.  Let $ P = \Phi \big( \overline{ M(\mathbf{p})} \big) $.  Anderson-Kogan \cite{jaredmisha} gave a combinatorial algorithm, called \textbf{collapse}, for use in calculating the vertices of $ P $ from the Kostant picture $ \mathbf{p} $.  

On the other hand, the above considerations show that $ P = P(M_\bigdot) $.  The values of $ M_\gamma $ for all $ \gamma \in \Gamma^\wi $ are linearly determined from $ n_\bigdot $.  The other values of $ M_\bigdot $ are determined by the tropical Pl\"ucker relations.  The positions of the vertices of $ P $ are determined from $ M_\bigdot $ by the usual vertex/hyperplane correspondence (\ref{eq:mufromM}).

Thus we have two combinatorial procedures for calculating the vertices $ \mu_\bigdot $ of $P $ from the $ \wi$-Lusztig datum $ n_\bigdot$: the collapse algorithm and the method of solving the tropical Pl\"ucker relations.  The Anderson-Kogan method is more explicit and perhaps easier to work with.  Both procedures produce the same answer, but in fact more is true --- we can understand the collapse algorithm as a series of applications of the tropical Pl\"ucker relations.  The remainder of this section will be devoted to the explanation of this statement.

Actually, the vertices produced by the Anderson-Kogan method and those which we produce differ in the labelling of the vertices by the Weyl group.  If $ \nu_\bigdot $ denotes the Anderson-Kogan vertices, then $ \nu_w = \mu_{w w_0} $.

Collapse along $ k $ takes a Kostant picture $ \mathbf{p} $ and produces another Kostant picture $ \mathbf{p'} $, whose ``loops'' are naturally labelled $ (a,b) $ with $ a < b $ and $ a \ne k , b \ne k $ (see \cite[section 2.4] {jaredmisha}).  

Anderson-Kogan defined collapse by a combinatorial algorithm.  As the definition is quite involved, we will not give it here.  However we can summarize the algorithm by the following algebraic statement.

\begin{Lemma} \label{th:ptop'}
If $ a < k < b $, then
\begin{equation} \label{eq:ptop'} 
p'_{(a,b)} = \min \Big( \sum_{r = k}^{b-1} p_{(a,r)} - p'_{(a,r)}, \sum_{s = a+1}^k p_{(s,b)} - p'_{(s,b)} \Big),
\end{equation}
where by convention $ p'_{(a,k)} = 0 = p'_{(k,b)} $.
If $ k < a $ or $ b < k $, then $ p'_{(a,b)} = p_{(a,b)} $.
\end{Lemma}

\begin{proof}
We give a sketch of the proof which will be comprehensible only to those familiar with the collapse algorithm.  Note that every loop with left edge at $ a $ and right edge at $ b$ (from now on called an $(a,b)$-loop) is produced as the join of an $ (a,r) $ loop and an $(s,b) $ for some $ k \le r < b $ and $ a < s \le k $.  Now every such $ (a,r) $ and $(s,b) $ loop is joined at some stage of the collapse algorithm, so the two parts of the $ \min $ represent the amount of $ (a,r) $ and $(s,b) $ loops not used to make smaller loops.  The production of $(a,b)$ loops by joining is then given by the minimum number of the available raw materials.
\end{proof}

In \cite[section 2.4]{jaredmisha}, collapse along $ k $ is used to understand the vertices $ \nu_w $ for all $ w $ such that $ w(1) = k $.  For us, these will be the vertices $ \mu_w $ with $ w(n) = k $.

The set $ \{ w \in W : w(n) = k \} $ forms a facet of the permutahedron.  This facet is naturally isomorphic to the permutahedron of $ SL_{n-1} $ and we consider a path along this facet corresponding to the above reduced word for $ SL_{n-1} $.  We extend this path to a reduced word $\wi' $ for $ SL_n $.  This reduced word gives us a labelling of the edges of the path by positive roots of $ SL_n $.  The labelling of the edges of the path lying in the facet is independent of how we extend the path outside of this facet.  In fact, they are always labelled by the positive roots $ (a,b) $ such that $ a \ne k, b \ne k $.

\begin{Lemma}
Let $ n_\bigdot $ be a reduced word for $ \wi $ and let $ n'_\bigdot $ be the corresponding reduced word for $ \wi' $.  For any positive root $ (a,b) $, let $ p_{(a,b)} = n_k $ where $ k $ is such that $ \beta_k^\wi = (a,b) $ and for any positive root $ (a,b) $ with $ a \ne k, b \ne k $, let $ p'_{(a,b)} = n'_k $ where $ k $ is such that $ \beta_k^{\wi'} = (a,b) $.  Then $ p'_\bigdot $ and $ p_\bigdot $ are related as in Lemma \ref{th:ptop'}.
\end{Lemma}

\begin{proof}
Assume $ a < k < b $.
Let $ M_\bigdot $ be the corresponding BZ datum.  By applying the usual conversion (\ref{eq:Mfromn}), we see that \begin{gather}
p_{(a,b)} = - M_{[a+1 \cdots b]} - M_{[a \cdots b-1]} + M_{[a+1 \cdots b-1]} + M_{[a \cdots b]} \label{eq:p} \\
p'_{(a,b)} = - M_{[a+1 \cdots \hat{k} \cdots b]} - M_{[a \cdots \hat{k} \cdots b-1]} + M_{[a+1 \cdots \hat{k} \cdots b-1]} + M_{[a \cdots \hat{k} \cdots b]} \label{eq:p'}.
\end{gather}

We expand out the sum in the RHS of (\ref{eq:ptop'}) using (\ref{eq:p}) and (\ref{eq:p'}).  Then we substitute (\ref{eq:p'}) into the LHS of (\ref{eq:ptop'}).  After cancelling some terms, we see that we must prove that for all $ a < k < b $, 
\begin{equation} \label{eq:tprcoll}
M_{[a \cdots \hat{k} \cdots b]} + M_{[a+1 \cdots b-1]} = \min \big( M_{[a+1 \cdots \hat{k} \cdots b]} + M_{[a \cdots b-1]}, M_{[a \cdots \hat{k} \cdots b-1]} + M_{[a+1 \cdots b]} \big).
\end{equation}
But this is exactly a tropical Pl\"ucker relation and so the result follows.
\end{proof}

Thus our theory gives the following interpretation of collapse.
\begin{Theorem} \label{th:mycollapse}
The Kostant picture produced by collapse along $ k $ gives the lengths of the edges along the above path inside the ``$ w(n) = k $'' facet of the MV polytope.  Moreover, the algorithm of collapse along $k $ is equivalent to recursively solving all of the tropical Pl\"ucker relations of the form (\ref{eq:tprcoll}).
\end{Theorem}
 
\subsection{Vertices}
Now that we have managed to see how the collapse algorithm fits into our setup, it remains only to understand the inductive way that Anderson-Kogan compute the vertices $ \nu_\bigdot $.  Fix $ v \in W $ and let $ \nu_v = (\nu^1, \dots, \nu^n) $.  First they compute $ \nu^k $ where $ k = v(1) $, and then they compute the rest of the components of the vertex using the Kostant datum produced by collapse along $ k $.  Since we understand collapse along $ k $, it remains only to understand their formula for $ \nu_k $. 

They show (see the proof of Theorem E in \cite{jaredmisha}) that for any $ L \in M^\ddag(\mathbf{p}) $, 
\begin{equation} \label{eq:theirvert}
\nu^k = \delta_k(L) + \dim_0(L) - \dim_0(L \cap U_\gamma),
\end{equation} 
where $ \gamma = [1 \cdots \hat{k} \cdots n]$.

This vertex $ \nu_v $ will be our vertex $ \mu_w $ where $ w = v w_0 $.  So $ w (n) = k $ which implies that $ w \cdot \fund_{n-1} = \gamma $.  Hence we have that
\begin{equation} \label{eq:ourvert}
\langle \mu_w , \gamma \rangle = M_\gamma.
\end{equation}
 
But if $ L \in A(M_\bigdot) $, then $ M_\gamma = \D_\gamma(L) = \rdim(L \cap U_\gamma) $ by Proposition \ref{th:dgam}.  Now we apply Lemma \ref{th:dim0} to conclude that
\begin{align*}
\rdim(L \cap U_\gamma) &= \dim_0(L \cap U_\gamma) + \sum_{i \ne k} \delta_i(L), \\
\rdim(L) &= \dim_0(L) + \sum_i \delta_i(L).
\end{align*}
Taking the difference of these two equations and combining with (\ref{eq:ourvert}) yields 
\begin{equation*}
\langle \mu_w, \gamma \rangle = \dim_0(L \cap U_\gamma) - \dim_0(L) - \delta_k(L)
\end{equation*}
which is equivalent to (\ref{eq:theirvert}) since $ \langle \nu, \gamma \rangle = - \nu^k $ as we are in the coweight lattice of $ SL_n$.

\appendix
\section{Pseudo-Weyl polytopes} \label{se:appendix}
The purpose of this appendix is to prove Proposition \ref{th:BZtoGGMS}.  To do so we will introduce the notion of dual fan to a polytope.  This will also put the concept of pseudo-Weyl polytope on a firmer footing.  

We thank A. Knutson, D. Speyer, and B. Sturmfels for explaining some of the concepts presented here and for suggesting this method of proving Proposition \ref{th:BZtoGGMS}.  Many of the definitions presented here can be found in \cite{E} and \cite{Z}.  The general results presented here (Theorems \ref{th:supp} and \ref{th:vert}) are known to experts but we could not find them in the literature.  A version of Theorem \ref{th:supp} appears in \cite[section 3.1]{F} in the context of ample line bundles on toric varieties.

Let $ V $ denote a real vector space and $ V^\star $ its dual.  We will be interested in the case $ V = \realt $.

\subsection{Support functions}
If $ P $ is a convex subset of $ V $,  we define the \textbf{support function} of P, $ \psi_P : V^\star \rightarrow \R $, by 
\begin{equation*}
\psi_P(\alpha) := \min_{v \in P} \langle v, \alpha \rangle.
\end{equation*}

The support function is a homogeneous, concave function on $ V^\star$, i.e. 
\begin{gather*}
\psi_P(\lambda \alpha) = \lambda \psi_P(\alpha) \text{ if } \lambda \in \R^\times \text{ and } \\
\psi_P(\alpha + \beta) \ge \psi_P(\alpha) + \psi_P(\beta).
\end{gather*}

Conversely, given any homogeneous, concave function $\psi $ on $ V^\star $, we can define the set 
\begin{equation*}
P(\psi) := \{ v \in V : \langle v, \alpha \rangle \ge \psi(v) \rangle \}.
\end{equation*}

The general theory of convexity tells us that these two maps are inverse to each other and set up a bijection
\begin{equation*}
\left \lgroup \begin{array}{c} \text{convex subsets} \\ \text{ of } $ V $ \end{array} \right \rgroup \longleftrightarrow 
\left \lgroup \begin{array}{c} \text{ homogeneous concave } \\ \text{ functions on }  V^\star  \end{array} \right \rgroup.
\end{equation*}

We will now proceed to examine a special case of this bijection.

\subsection{Fans}
A \textbf{polyhedral cone} in $ V^\star $ is a finite intersection of closed linear half spaces.

A (complete) \textbf{fan} $\mathcal{F}$ in $V^\star$ is a finite collection of nonempty polyhedral cones of $ V^\star $ such that 
\begin{enumerate}
\item every nonempty face of a cone in $\mathcal{F}$ is also a cone in $ \mathcal{F} $,
\item the intersection of any two cones in $ \mathcal{F} $ is a face of both, and
\item the union of all the cones in $ \mathcal{F} $ is $ V^\star$.
\end{enumerate}

A fan $ \mathcal{F} $ induces an equivalence relation on $ V^\star $ whose equivalence classes are the interiors of the cones of $ \mathcal{F} $.  The fan can be recovered from this equivalence relation, thus we can view fans as a special class of equivalence relation on $ V^\star $. 

We will mostly be concerned with the \textbf{Weyl fan} in $ \realt^\star $.  The maximal cones of this fan are the cones
\begin{equation*}
C_w^\star := \{ \alpha \in \realt^\star : \langle w \cdot \alpha_i^\vee, \alpha \rangle \ge 0 \text{ for all } i \}.
\end{equation*}
All other cones are obtained by intersecting these maximal cones.  

Given a polytope $ P $ in $ V $ (for us polytopes are always assumed to be convex),
we can construct a fan in $ V^\star $ called the dual fan $ \mathcal{N}(P) $ of $ P $.  For $ \alpha \in V^\star$, let  
\begin{equation*}
 M(P, \alpha) = \{ v \in P : \langle v , \alpha \rangle = \psi_P(\alpha) \}
\end{equation*}
be the subset of $ P $ where $ \langle \cdot, \alpha \rangle $ is minimized.  Note that this subset will always be a face of $ P $.
The cones $C^\star_F $ of $ \mathcal{N}(P) $ are indexed by the faces $ F $ of $ P $ and are given by
\begin{equation*}
C_F^\star := \{ \alpha \in V^\star : F \subset M(P, \alpha) \}.
\end{equation*}
The corresponding equivalence relation is 
\begin{equation*}
\alpha \sim \beta \Leftrightarrow M(P, \alpha) = M(P, \beta)
\end{equation*}

\begin{Proposition}
The dual fan of the permutahedron is the Weyl fan.
\end{Proposition}

\begin{proof}
Recall that for each $ w \in W $, $ w w_0 \cdot \rho^\vee $ is a vertex of the permutahedron.  We will show that the cone $ C_{w w_0\cdot \rho^\vee}^\star $ dual to this vertex is $ C_w^\star $.  

Since the local cone of the permutahedron at the $ w w_0\cdot \rho^\vee $ vertex is $ C_w^{w w_0\cdot \rho^\vee} $, we have that 
\begin{equation*}
C_{w w_0 \cdot \rho^\vee}^\star = \{ \alpha \in \realt^\star : \langle v, \alpha \rangle \ge \langle w w_0\cdot \rho^\vee, \alpha \rangle \text{ for all } v \in C_w^{w w_0\cdot \rho^\vee} \}
\end{equation*}

So $ \alpha $ lies in the dual cone if and only if $ \langle v, \alpha \rangle \ge 0  $ for all $ v \in C_w^0$.  Now the cone $ C_w^0 $ is spanned by $ w \cdot \fund_i $ for all $ i \in I $ and so $ \alpha $ is in the dual cone if and only if $ \langle w \cdot \alpha_i^\vee, \alpha \rangle \ge 0 $ for all $ i \in I $ as desired.
\end{proof}

A fan $\mathcal{F}_1 $ is said to be a \textbf{coarsening} of a fan $ \mathcal{F}_2 $ if every cone of $ \mathcal{F}_1 $ is a union of cones of $ \mathcal{F}_2 $.  Equivalently, the equivalence relation $ \underset{1}{\sim} $ corresponding to $ \mathcal{F}_1 $ is stronger than the equivalence relation $ \underset{2}{\sim} $ corresponding to $ \mathcal{F}_2$, i.e. $ \alpha \underset{2}{\sim} \beta \Rightarrow \alpha \underset{1}{\sim} \beta$.

A polytope $P$ is called an $\mathcal{F}$-\textbf{polytope} if its dual fan is a coarsening of $\mathcal{F}$.

With these notions in hand, we can now give a better definition of pseudo-Weyl polytope.  A \textbf{pseudo-Weyl polytope} is a polytope in $ \realt $ with vertices in $ \cwl $, whose dual fan is a coarsening of the Weyl fan.  Later we will show that this definition agrees with our old one.

\subsection{Support functions of $\mathcal{F}$-polytopes}
We would like to see how to characterize $ \mathcal{F}$-polytopes in terms of their support functions.

A homogeneous, concave function $ \psi $ is said to be \textbf{linear} on $ \mathcal{F} $ if, whenever $ \alpha \sim \beta $, we have $\psi(\alpha) + \psi(\beta) = \psi(\alpha + \beta) $.  Since a concave function is automatically continuous, this implies that the restriction of $\psi $ to any cone of $\mathcal{F} $ is linear.

\begin{Theorem} \label{th:supp}
The maps $ P \mapsto \psi_P $ and $ \psi \mapsto P(\psi) $ give a bijection
\begin{equation*}
\left \lgroup \mathcal{F}-\text{Polytopes} \right \rgroup \longleftrightarrow 
\left \lgroup \begin{array}{c} \text{homogeneous, concave functions} \\ \text{ which are linear on $\mathcal{F}$} \end{array} \right \rgroup.
\end{equation*}
\end{Theorem}

\begin{proof}
Since these maps are inverses to each other we just need to check that if $ P $ is an $\mathcal{F}$-polytope, then $ \psi_P $ is linear on $ \mathcal{F} $ and conversely if $ \psi $ is a homogeneous, concave function, linear on $ \mathcal{F}$, then $ P(\psi) $ is an $\mathcal{F}$-polytope.

First, assume that $ P $ is an $ \mathcal{F}$-polytope.  Let $\alpha \sim \beta $.  Then $ \alpha $ and $ \beta $ are also equivalent under the $\mathcal{N}(P) $ equivalence relation (since $ \mathcal{N}(P) $ is a coarsening of $ \mathcal{F} $).  So $ M(P, \alpha) = M(P, \beta) $.  Hence there exists $ v \in P $ such that $ \langle v, \alpha \rangle = \psi_P(\alpha) $ and $ \langle v, \beta \rangle = \psi_P(\beta) $.  

Hence, $ \psi_P(\alpha + \beta) \le \langle v, \alpha + \beta \rangle = \psi_P(\alpha) + \psi_P(\beta)$.  Hence $ \psi_P(\alpha + \beta) = \psi_P(\alpha) + \psi_P(\beta) $ as desired.  So $ \psi_P $ is linear on $ \mathcal{F} $.

Now assume that $ \psi$ is a homogeneous, concave function which is linear on $ \mathcal{F}$.  Let $ \alpha \sim \beta $ in $ \mathcal{F} $.  We would like to show that $ M(P, \alpha) = M(P, \beta) $ since this will show that $ \alpha $ and $ \beta $ are similar under the $ \mathcal{N}(P) $ equivalence relation.

Suppose that $ v \in M(P, \alpha) $ but $ v \notin M(P, \beta) $, so $\psi(\beta) < \langle v, \beta \rangle $.  Since the equivalence classes of $ \mathcal{F} $ are the interiors of cones, there exists $ t > 0 $ such that $ \alpha - t \beta \sim \beta $.  By linearity $ \psi(\alpha - t \beta) + \psi( t \beta) = \psi(\alpha) $. Hence,
\begin{equation*}
\langle v, \alpha - t \beta \rangle + \langle v, t \beta \rangle > \psi(\alpha - t \beta) + t \psi(\beta) = \psi(\alpha) = \langle v, \alpha \rangle
\end{equation*}
which is a contradiction.

So we conclude that $ M(P, \alpha) \subset M(P, \beta) $ and similarly $ M(P, \beta) \subset M(P, \alpha) $.  Hence $ \alpha $ and $ \beta $ are similar under the $ \mathcal{N}(P) $ equivalence relation.
\end{proof}

\subsection{Vertex data}
If $ P $ is a polytope, then there is a natural bijection between the vertices of $ P $ and the maximal cones of $ \mathcal{N}(P) $.

Let $ \mathcal{F} $ be a fan and let $ P $ be an $\mathcal{F}$-polytope.  Let $ \{ C^\star_x : x \in X \} $ be the set of maximal cones of $\mathcal{F}$.  Each maximal cone of $ \mathcal{F} $ is a subcone of a unique maximal cone of $ \mathcal{N}(P) $ and so we get a surjective map
\begin{equation*}
X \twoheadrightarrow \text{max cones of } \mathcal{N}(P) = \text{vertices of } P.
\end{equation*}
Let $ p_x $ denote the image of $ x \in X $ under this map.  

The collection $ \Verti(P) := \big( p_x \big)_{x \in X} $ is called the \textbf{vertex data} of $P $.

For each $ x \in X $ define a partial order $ \le_x $ on $ V $ by 
\begin{equation*}
v \le_x w \Leftrightarrow \langle v, \alpha \rangle \le \langle w, \alpha \rangle \text{ for all } \alpha \in C_x^\star
\end{equation*}

Suppose we have a collection of points of $ V $, $ p_\bigdot = (p_x)_{x \in X} $ such that $ p_y \ge_x p_x $ for all $ x,y \in X $.  Then we define
\begin{equation*}
P(p_\bigdot) := \{ v \in V : v \ge_x p_x \text{ for all } x \in X \}.
\end{equation*}

\begin{Theorem} \label{th:vert}
The maps $ P \mapsto \Verti(P) $ and $ p_\bigdot \mapsto P(p_\bigdot) $ give a bijection
\begin{equation*}
\left \lgroup \mathcal{F}-\text{polytopes} \right \rgroup \longleftrightarrow 
\left \lgroup \begin{array}{c} \text{ collections } (p_x)_{ x \in X } \text{ such that } \\
  p_y \ge_x p_x \text{ for all }  x,y \in X  \end{array} \right \rgroup
\end{equation*}

Moreover, if $ P $ and $p_\bigdot $ correspond under this bijection then the support function of $ P $ satisfies 
\begin{equation*}
\psi_P(\alpha) = \langle p_x, \alpha \rangle \text{ if }  \alpha \in C_x^\star.
\end{equation*}
\end{Theorem}

\begin{proof}
First, we would like to show that if $ P $ is an $\mathcal{F}$-polytope, then $ p_\bigdot := \Verti(P) $ satisfies the desired condition.  Let $ x,y \in X $ and let $ \alpha \in C_x^\star $.  By definition, $ p_x \in M(P, \alpha) $.  So
$ \psi_P(\alpha) = \langle p_x, \alpha \rangle $.  Since $ p_y \in P $,
\begin{equation*}
\langle p_y, \alpha \rangle \ge \psi_P(\alpha) = \langle \alpha, p_x \rangle.
\end{equation*}
Since this holds for all $ \alpha \in C_x^\star $, $ p_y \ge_x p_x $ as desired.

Now, we show that $ P(\Verti(P)) = P $.  Note that $ v \ge_x p_x $ if and only if $ \langle v, \alpha \rangle \ge \langle p_x, \alpha \rangle = \psi_P(\alpha) $ for all $ \alpha \in C_x^\star $.  So we see that 
\begin{equation*}
P(p_\bigdot) = \{ v \in V : \langle v, \alpha \rangle \ge \psi_P(\alpha) \text{ for all } \alpha \in V^\star \}.
\end{equation*}
So $ P(p_\bigdot) = P(\psi_P) = P $ as desired.

Next, we would like to show that if $ p_\bigdot $ satisfies the hypothesis, then $ P(p_\bigdot) $ is an $\mathcal{F}$-polytope and $\Verti(P(p_\bigdot)) = p_\bigdot $.  Define a function $\psi : V^\star \rightarrow \R $ by 
\begin{equation*}
\psi(\alpha) = \langle p_x, \alpha \rangle \text{ if }  \alpha \in C_x^\star . 
\end{equation*}
 To see that $ \psi $ is well-defined, suppose that $ \alpha \in C_x^\star $ and $ \alpha \in C_y^\star $.  Then since $ p_y \ge_x p_x $, $ \langle p_y, \alpha \rangle \ge \langle p_x, \alpha \rangle $.  Similarly the opposite inequality holds and so we see that $ \langle p_y, \alpha \rangle = \langle p_x, \alpha \rangle $.

Now, we claim that $ \psi $ is homogeneous and concave.  Homogeneity is clear.  Suppose that $ \alpha \in C_x^\star, \beta \in C_y^\star $ for some $ x, y \in X $.  Then there exists $ u \in X $ such that $ \alpha + \beta \in C_u^\star $. So,
\begin{equation*}
\psi(\alpha + \beta) = \langle p_u, \alpha + \beta \rangle = \langle p_u, \alpha \rangle + \langle p_u, \beta \rangle \ge \psi(\alpha) + \psi(\beta)
\end{equation*}
where the inequality follows from the fact that $ p_u \ge_x p_x $ and $p_u \ge_y p_y $.  So $ \psi $ is concave.  Finally, we claim that $\psi $ is linear on $ \mathcal{F}$.  Suppose that $ \alpha \sim \beta $.  Then there exists $x \in X $ such that $ \alpha, \beta \in C_x^\star $.  Then $ \alpha + \beta $ is also in $ C_x^\star $.  So $ \psi(\alpha + \beta) = \langle p_x, \alpha + \beta \rangle = \psi(\alpha) + \psi(\beta) $ as desired.

Hence $ P(\psi) $ is an $ \mathcal{F} $-polytope.  But we have already seen that $ P(\psi) = P(p_\bigdot) $ and so $ P(p_\bigdot) $ is an $ \mathcal{F} $-polytope.  Moreover, we already saw that if $ p'_\bigdot = \Verti(P(\psi)) $, then $ \psi(\alpha) = \langle p'_x, \alpha \rangle $ for all $ \alpha \in C_x^\star$.  Since the cone $C_x^\star $ is maximal, it spans $ V^\star $ and hence $ p'_x = p_x $ for all $ x $ as desired. 
\end{proof}

\begin{Corollary}
Our two definitions of pseudo-Weyl polytope agree.
\end{Corollary}

Note that we have also proven the first part of Proposition \ref{th:BZtoGGMS}.

\subsection{Hyperplane data}
A polyhedral cone $ C $ in $ V^\star $ is called \textbf{simplicial} if there exists a basis $ \alpha_1, \dots, \alpha_n $  for $ V^\star $ such that $ C = \{ \lambda_1 \alpha_1 + \dots + \lambda_n \alpha_n : \lambda_i \ge 0 \} $.  These $ \alpha_1, \dots, \alpha_n $ will necessarily be along the rays (one-dimensional faces) of the cone.
A fan $ \mathcal{F} $ is called \textbf{simplicial} if all of its cones are simplicial.  For example, the Weyl fan is simplicial since the cone $ C_w^\star $ is spanned by the vectors $ \{ w \cdot \fund_i : i \in I \} $.  

From now on, we assume that $ \mathcal{F} $ is simplicial and let $ \Gamma $ be a set of vectors, one lying in each ray of $ \mathcal{F} $.  So for any cone $ C \in \mathcal{F} $, $ C $ is the positive linear span of the vectors $ \Gamma \cap C $.  For example, when $ \mathcal{F} $ is the Weyl fan, the set of chamber weights $ \Gamma $ is such a set.  

If $ \psi $ is a homogeneous concave function, linear on $ \mathcal{F} $, then $ \psi $ is determined by its restriction to the rays of $ \mathcal{F} $.   Hence we get a sequence of real numbers $ \big( M_\gamma := \psi(\gamma) \big)_{\gamma \in \Gamma}$ which determine $ \psi $.  Moreover, in this case, we see that 
\begin{equation*}
P(\psi) = \{ v \in V : \langle v, \gamma  \rangle \le M_\gamma \text{ for all } \gamma \in \Gamma \}.
\end{equation*}
The collection $ M_\bigdot $ is called the \textbf{hyperplane datum} of $ P(\psi) $.

Conversely, given a sequence of real numbers $ \big( M_\gamma \big)_{\gamma \in \Gamma} $, we can ask if there exists a homogeneous, concave function $\psi $ which is linear on $ \mathcal{F} $ such that $ \psi(\gamma) = M_\gamma $ for all $ \gamma $.   In fact, such a sequence always defines a function $ \psi_{M_\bigdot} $ in the following way.  Since $ \mathcal{F} $ is simplicial, every $ \alpha \in V^\star $ can be written uniquely as a positive linear combination $ \alpha = \lambda_1 \gamma_1 + \dots + \lambda_n \gamma_n $  with $ \gamma_i \in \Gamma $ and $ \lambda_i \ge 0 $.  Then we define 
\begin{equation*}
\psi_{M_\bigdot}(\alpha) := \lambda_1 M_{\gamma_1} + \dots + \lambda_n M_{\gamma_n}.
\end{equation*}
Note that $ \psi_{M_\bigdot} $ is homogeneous and linear on $ \mathcal{F} $.  However, it will not always be true that $ \psi_{M_\bigdot} $ is concave.

\begin{Lemma} \label{th:nondeg}
If $ \mathcal{F} $ is the Weyl fan and $ \Gamma $ is the set of chamber weights, then $ \psi_{M_\bigdot} $ is concave if and only if $ M_\bigdot $ satisfies the edge inequalities (\ref{eq:nondeg}).
\end{Lemma}

\begin{proof}
First, we show that if $ \psi_{M_\bigdot} $ is concave, then $ M_\bigdot $ satisfies the edge inequalities.  For any $ i \in I $ and $ w \in W $, note that
\begin{equation*}
w s_i \cdot \fund_i + w \cdot \fund_i  =  \sum_{j \ne i} -a_{ji} w \cdot \fund_j.
\end{equation*}

Since $w \cdot \fund_j $ all lie in the same cone of the Weyl fan and since $ a_{ji} \le 0 $ for $ j\ne i$, by linearity, homogeneity, and concavity, we have that
\begin{equation*}
\psi( w s_i \cdot \fund_i) + \psi(w \cdot \fund_i) \le \sum_{j \ne i} - a_{ji} \psi(w \cdot \fund_j).
\end{equation*}
This implies that edge inequality among the $ M_\bigdot $.

Conversely, assume that $ M_\bigdot $ satisfies the edge inequalities, and define $ \psi_{M_\bigdot} $ as above.  We would like to show that $ \psi_{M_\bigdot} $ is concave.  

For each $ w \in W $, let $\psi_w $ denote the unique linear function on $ \realt^\star $ such that $\psi_w(w \cdot \fund_i) =  M_{w \cdot \fund_i} $ for all $ i \in I $.  So $ \psi $ and $ \psi_w $ agree on $ C_w^\star $.  By the same argument as in the second half of the proof of Theorem \ref{th:vert}, it suffices to show that $ \psi(\alpha) \le \psi_w(\alpha) $ for all $ \alpha \in \realt^\star $.  To prove this it suffices to show that $ \psi(\gamma) \le \psi_w(\gamma) $ for all $ \gamma \in \Gamma $ and $ w \in W$. 

For simplicity, we will prove this last statement for $ \gamma = \fund_k $ for some $ k $.  Our proof will proceed by induction on $ W $ using the weak Bruhat order.  To prove the statement for general $ \gamma $ requires a different partial order adapted to $ \gamma $.

The base case of $ w = e $ is clear.  So assume $ w \in W $, $ w \ne e $ and that $\psi(\gamma) \le \psi_v(\gamma) $ for all $ v < w $ in the weak Bruhat order.

Let $ \lambda_j \in \R $ be such that $ \gamma = \sum \lambda_j w \cdot \fund_j $.  Since $ w \ne e $, there exists $ i \in I $ such that $ w s_i < w $.  Hence $ w \cdot \alpha_i^\vee $ is a negative coroot.  So $ \langle w \cdot \alpha_i^\vee, \fund_k \rangle \le 0 $ and hence $ \lambda_i \le 0 $.

Since $ \fund_i = - s_i \cdot \fund_i - \sum_{j \ne i} a_{ji} \fund_j $ and $ s_i \cdot \fund_j = \fund_j $ for $ j \ne i $, we see that 
\begin{equation*}
\gamma = \sum_{j \ne i} (\lambda_j - \lambda_i a_{ji} ) w s_i \cdot \fund_j - \lambda_i w s_i \cdot \fund_i.
\end{equation*}

Since $ w s_i < w $, by induction we have that $ \psi(\gamma) \le \psi_{w s_i}(\gamma) $ and so
\begin{equation} \label{eq:expandedM}
M_\gamma \le \sum_{j \ne i} (\lambda_j - \lambda_i a_{ji} ) M_{w s_i \cdot \fund_j} - \lambda_i M_{w s_i \cdot \fund_i} = \sum_{j \ne i} \lambda_j M_{w \cdot \fund_j} - \lambda_i(M_{w s_i \cdot \fund_i} + \sum_{j \ne i} a_{ji} M_{w \cdot \fund_j} ).
\end{equation}

Now the edge inequality tells us that 
\begin{equation*} 
M_{w s_i \cdot \fund_i} +  \sum_{j \ne i} a_{ji} M_{w \cdot \fund_j} \le - M_{w \cdot \fund_i} 
\end{equation*}

So multiplying this equation by $ -\lambda_i $ and combining with (\ref{eq:expandedM}) shows that 
\begin{equation*}
M_\gamma \le \sum_j \lambda_j M_{w \cdot \fund_j}
\end{equation*}
as desired.

Hence we have proven the statement for $ w $.  This completes the induction argument.
\end{proof}

Let $ P $ be a polytope whose dual fan is a coarsening of the Weyl fan.  Let $\mu_\bigdot $ be its vertex data and $ M_\bigdot $ be its hyperplane data.  Then by Theorem \ref{th:vert}, they are related by
\begin{equation*}
M_{w \cdot \fund_i} = \langle \mu_w , w \cdot \fund_i \rangle.
\end{equation*}
So we see that $ M_\gamma \in \Z $ for all $ \gamma $ if and only if $ \mu_w \in \cwl $ for all $ w $.

Combining Theorem \ref{th:supp}, Theorem \ref{th:vert}, Lemma \ref{th:nondeg} and the above remark, gives the proof of Proposition \ref{th:BZtoGGMS}.

\subsection{Minkowski sums of pseudo-Weyl polytopes}
We close this section with the proof of Lemma \ref{th:minksum} concerning Minkowski sums of pseudo-Weyl polytopes.

\begin{proof}[Proof of Lemma \ref{th:minksum}]
If $ A, B $ are polytopes, then the dual fan of the Minkowski sum $ A + B $ is the common refinement of the two dual fans $ \mathcal{N}(A), \mathcal{N}(B) $ (see \cite[Prop 7.12]{Z}).  If two fans are both coarsenings of the Weyl fan, then so is their common refinement.  Hence the Minkowski sum of pseudo-Weyl polytopes is again a pseudo-Weyl polytope.  Moreover, it is clear that the support function of $ A+B$ is $ \psi_A + \psi_B $ and so the second half of the result follows.
\end{proof}

\end{document}